\documentclass[11pt]{article}

\usepackage[margin=1.15in]{geometry}
\usepackage{amsmath,amssymb,amsthm,mathtools}
\usepackage[shortlabels]{enumitem}
\usepackage{microtype}
\usepackage[colorlinks=true,linkcolor=blue,citecolor=blue,urlcolor=blue]{hyperref}
\usepackage[capitalise,noabbrev]{cleveref}
\usepackage{tikz-cd}

\newtheorem{theorem}{Theorem}[section]

\theoremstyle{definition}
\newtheorem{definition}[theorem]{Definition}
\newtheorem{definition-and-lemma}[theorem]{Definition and Lemma}

\theoremstyle{remark}
\newtheorem{remark}[theorem]{Remark}

\crefname{enumi}{item}{items}
\Crefname{enumi}{Item}{Items}
\crefname{equation}{}{}
\Crefname{equation}{}{}

\ExplSyntaxOn

\NewDocumentEnvironment {athm} {m m o} {
\str_if_eq:noTF {example} {#1} {
  \bool_gset_true:N \g_example_bool
} {
  \bool_gset_false:N \g_example_bool
}
\IfNoValueTF{#3}{
\begin{#1}\label{#2}\global\def\loc{#2}
}{
\begin{#1}[#3]\label{#2}\global\def\loc{#2}
}
}{
\end{#1}
}

\NewDocumentEnvironment {adef} {m} {
\begin{definition}\label{#1}\global\def\loc{#1}
}{
\end{definition}
}

\NewDocumentEnvironment{aproof} { s } {
\bool_if:NTF \g_example_bool {
  \begin{proof}[Proof~for~\cref{\loc}]
} {
  \begin{proof}[Proof~of~\cref{\loc}]
}
\bool_gset_false:N \g_finishproof_bool
}{
  \IfBooleanT{#1}{
\bool_if:NTF \g_finishproof_bool {}
{\finishproofthus}}
\end{proof}
}

\NewDocumentCommand{\finishproofthus} {} {
  \bool_gset_true:N \g_finishproof_bool 
  \bool_if:NTF \g_example_bool {
    The~proof~for~\cref{\loc}~is~thus~complete.
  } {
    The~proof~of~\cref{\loc}~is~thus~complete.
  }
}
\NewDocumentCommand{\finishproofthis} {} {
  \bool_gset_true:N \g_finishproof_bool 
  \bool_if:NTF \g_example_bool {
    This~completes~the~proof~for~\cref{\loc}.
  } {
    This~completes~the~proof~of~\cref{\loc}.
  }
}

\ExplSyntaxOff

\newcommand{\llabel}[1]{\label{\loc.#1}}

\ExplSyntaxOn

\clist_new:N \l_localreflist
\clist_new:N \l_reflist

\NewDocumentCommand{\lref} { s m } {
  \clist_set:No \l_localreflist {#2}
  \clist_clear:N \l_reflist
  \clist_map_inline:Nn \l_localreflist { \clist_put_right:Nn \l_reflist {\loc.##1} }
  \IfBooleanTF {#1} {
    \labelcref{\l_reflist}
  } {
    \cref{\l_reflist}
  }
}

\NewDocumentCommand{\Lref} { m } {
  \clist_set:No \l_localreflist {#1}
  \clist_clear:N \l_reflist
  \clist_map_inline:Nn \l_localreflist { \clist_put_right:Nn \l_reflist {\loc.##1} }
    \Cref{\l_reflist}
}

\NewDocumentCommand{\itref}{ s m m }{
  \clist_set:No \l_localreflist {#3}
  \clist_clear:N \l_reflist
  \clist_map_inline:Nn \l_localreflist { \clist_put_right:Nn \l_reflist {#2.##1} }
  \IfBooleanTF {#1} {
    \labelcref{\l_reflist}~in~\cref{#2}
  } {
    \cref{\l_reflist}~in~\cref{#2}
  }
}

\NewDocumentCommand{\Itref}{ m m }{
  \clist_set:No \l_localreflist {#2}
  \clist_clear:N \l_reflist
  \clist_map_inline:Nn \l_localreflist { \clist_put_right:Nn \l_reflist {#1.##1} }
  \Cref{\l_reflist}~in~\cref{#1}
}

\ExplSyntaxOff

\DeclareMathOperator{\im}{im}
\DeclarePairedDelimiter{\pr}{\lparen}{\rparen}
\DeclarePairedDelimiter{\br}{\lbrack}{\rbrack}
\providecommand{\given}{\colon}
\DeclarePairedDelimiterX{\set}[1]{\lbrace}{\rbrace}{%
  #1%
}

\DeclarePairedDelimiter{\ang}{\langle}{\rangle}

\DeclareMathOperator{\Tor}{Tor}

\NewDocumentCommand{\shortsum}{e{_^}}{%
  \mathchoice
    {\shortsumaux{\displaystyle}{#1}{#2}}
    {\shortsumaux{\textstyle}{#1}{#2}}
    {\shortsumaux{\scriptstyle}{#1}{#2}}
    {\shortsumaux{\scriptscriptstyle}{#1}{#2}}%
}
\makeatletter
\newcommand{\shortsumaux}[3]{%
  \mathop{\raisebox{0pt}[.9\height][.9\depth]{%
    $#1\m@th\sum
      \IfValueT{#2}{_{#2}}%
      \IfValueT{#3}{^{#3}}$%
  }}%
}
\makeatother

\newcommand{\Z}{\mathbb Z}
\newcommand{\Q}{\mathbb Q}
\newcommand{\R}{\mathbb R}
\newcommand{\N}{\mathbb N}

\newcommand{\NS}[2]{\widehat{\Z #1}_{#2}}
\newcommand{\G}{\Gamma}
\newcommand{\FP}{\mathrm{FP}}
\newcommand{\F}{\mathrm F}
\newcommand{\id}{\mathrm{id}}

\newcommand{\tensor}{\otimes}

\newcommand{\lhdnormal}{\mathrel{\unlhd}}

\newcommand{\Nov}[2]{\Sigma_{\mathrm{Nov}}^{#1}(#2;\Z)}
\newcommand{\Sig}[2]{\Sigma^{#1}(#2;\Z)}

\newcommand{\ses}[2]{\begin{gathered}#1\\#2\end{gathered}}

\title{Finiteness Properties of Fibre Products\\
over Virtually Nilpotent Quotients}
\author{Benno Kuckuck\\[0.5em]
\small Applied Mathematics: Institute for Analysis and Numerics,\\\small University of Münster, Germany; email: \url{bkuckuck@uni-muenster.de}}
\date{8 September 2026}

\begin{document}
\maketitle

\begin{abstract}
The $n$-$(n+1)$-$(n+2)$ theorem, recently established by Cohen and
Shusterman, says that if two groups of type $\F_{n+1}$ map onto a common
quotient $Q$ of type $\F_{n+2}$ and one of the two kernels is of type
$\F_n$, then the associated fibre product is of type $\F_{n+1}$.  We prove a
stronger, symmetric variant of this theorem when $Q$ is virtually nilpotent.
More precisely, for $i\in\{1,2\}$, let
$
  1\to N_i\to \G_i\xrightarrow{\pi_i}Q\to1
$
be short exact sequences of groups and let $P=\{(\gamma_1,\gamma_2)\in\Gamma_1\times\Gamma_2
\colon \pi_1(\gamma_1)=\pi_2(\gamma_2)\}$ be their fibre
product. Let $k,l,m,n\in\N_0$, assume that $N_1$ and $N_2$
are of type $\FP_k$ and $\FP_l$, respectively, and that
$\Gamma_1$ and $\Gamma_2$ are of type $\F_m$ and $\F_n$,
respectively. We prove that $P$ is of type $\F_{\min\{k+l+1,m,n\}}$.
The same formula holds with the homological finiteness properties $\FP_r$ in
place of $\F_r$.
This provides an alternative
proof of the homotopical and homological Virtual Surjections Theorem,
which gives a finiteness criterion for subgroups of direct products
in terms of their embedding in the ambient product.
\end{abstract}

\section{Introduction}

Taking direct products is one of the simplest ways of constructing new groups from old ones, but it has long been observed that, when the factors are ``far from abelian'', the subgroups of such direct products can exhibit surprisingly wild behaviour. Already in 1958, Mihaĭlova~\cite{mihailova1958entrance} showed that a direct product of two free groups of rank $2$ contains a finitely generated subgroup with undecidable membership problem. In 1971, Miller~\cite{miller1971group} showed that such direct products contain finitely generated subgroups with undecidable conjugacy problem. Remarkably, these phenomena do not happen for \emph{finitely presented} subgroups of such direct products: Baumslag and Roseblade~\cite{BaumslagRoseblade1984} showed that finitely presented subgroups of a direct product of two finitely generated free groups are themselves virtually direct products of two finitely generated free groups and are moreover quite rigidly embedded in the ambient product.

In parallel, Stallings~\cite{Stallings1963} constructed a finitely presented group that fails to have finitely generated integral homology in dimension $3$ and Bieri~\cite{Bieri1976} identified this as the first stage in a systematic family that arises in direct products of three or more free groups: A direct product of $n$ free groups of rank $2$ contains a subgroup that has strong finiteness exactly up to dimension $n-1$, but not above. More precisely, these groups are of type $\F_{n-1}$ but not of type $\FP_{n}$.

This suggested that higher finiteness properties may serve as a useful measure of complexity of subgroups of direct products of free groups. In a landmark series of papers, Bridson, Howie, Miller, and Short \cite{BHMS2002,BHMS2009,BridsonMiller2009,BHMS2013} developed this viewpoint with remarkable force: For direct products of free groups, surface groups, and, ultimately, limit groups, they established strong connections between higher finiteness properties and subgroup structure. Notably, they showed that every subgroup in a direct product of $n$ limit groups that is of type $\FP_n(\Q)$ is virtually itself a direct product of up to $n$ limit groups, settling a question of Sela. For finitely presented subgroups they obtained a much more refined description, leading to a characterization of finitely presented residually free groups and to solutions of several associated algorithmic problems.

A key ingredient of that proof was the observation that for a subgroup $G\leq \Gamma_1\times\Gamma_2\times\dots\times\Gamma_n$ of a direct product of $n\geq 2$ arbitrary finitely presented groups, finite presentability is already forced by a certain embedding property: If for every pair $i,j\in\{1,2,\dots,n\}$ of different indices, the image of the projection map $G\to\Gamma_i\times\Gamma_j$ has finite index in $\Gamma_i\times \Gamma_j$ ($G$ ``virtually surjects to pairs of factors''), then $G$ is finitely presented~\cite[Theorem~A]{BHMS2013}. This follows from a foundational
statement about finite presentability of fibre products, the $1$-$2$-$3$ theorem~\cite[Theorem~B]{BHMS2013}, which says that for a pair of short exact sequences of groups
\begin{equation}
  \label{eq:seq}
  \ses{
  1\to N_1\to\G_1\xrightarrow{\pi_1}Q\to1
  } {
  1\to N_2\to\G_2\xrightarrow{\pi_2}Q\to1
  }
\end{equation}
where $N_1$ is finitely generated, $\G_1$ and $\G_2$ are finitely presented, and $Q$ is of type $\F_3$,
the associated fibre product
\[
  P=\{(\gamma_1,\gamma_2)\in\G_1\times\G_2
       \colon \pi_1(\gamma_1)=\pi_2(\gamma_2)\}
\]
is finitely presented.
A natural question, asked in \cite{Kuckuck2014}, is whether this
phenomenon extends to higher finiteness properties. It was finally answered positively by
Cohen and Shusterman~\cite{CohenShusterman2026}.
Their $n$-$(n+1)$-$(n+2)$ theorem \cite[Corollary~1.8]{CohenShusterman2026} says that if, in \cref{eq:seq}, 
$N_1$ is of type $\F_{n}$, $\G_1$ and $\G_2$ are of type $\F_{n+1}$ and 
$Q$ is of type $\F_{n+2}$ then the
associated fibre product is of type $\F_{n+1}$.

\paragraph{Results.}
In the case where $Q$ is abelian, \cite{Kuckuck2014} proved a stronger, symmetric conclusion:
If, in \cref{eq:seq}, $N_1$ is of type $\F_k$,
$N_2$ is of type $\F_l$ with $k+l\geq n$, $\G_1$ and $\G_2$ are of type $\F_{n+1}$, and $Q$ is 
(virtually) finitely generated abelian, then the
associated fibre product is of type $\F_{n+1}$. Hence, in the
abelian quotient case, \emph{both} kernels may
contribute to the finiteness length of the fibre product.
The purpose of this paper is to extend this symmetric phenomenon from
virtually abelian to virtually nilpotent quotients.
\begin{athm}{theorem}{thm:main-fp}[Symmetric $n$-$(n+1)$-$(n+2)$ theorem for
  fibre products over
  virtually nilpotent quotients]
Let $m,n,k,l\in\N_0$,
let
\[
\ses{
    1\to N_1\to\G_1\xrightarrow{\pi_1}Q\to1
} {
    1\to N_2\to\G_2\xrightarrow{\pi_2}Q\to1
}
\]
be short exact sequences of groups, assume that
$N_1$ is of type $\FP_k$, that $N_2$ is of type $\FP_l$,
that $\G_1$ is of type $\FP_m$ (resp.\ $\F_m$), that $\G_2$ is of type $\FP_n$ (resp.\ $\F_n$), 
and that $Q$ is virtually nilpotent.  Then the fibre product
\[
  P=\{(\gamma_1,\gamma_2)\in\G_1\times\G_2
       \colon \pi_1(\gamma_1)=\pi_2(\gamma_2)\}
\]
is of type $\FP_{\min\{m,n,k+l+1\}}$ (resp.\ $\F_{\min\{m,n,k+l+1\}}$).
\end{athm}

This nilpotent case of the general $n$-$(n+1)$-$(n+2)$ theorem has particular significance,
as it was observed in \cite[Theorem~3.10]{Kuckuck2014} in the homotopical setting
and by Kochloukova and Lima~\cite[Theorem~F]{KochloukovaLima2018} in the homological setting that it suffices for establishing the higher-dimensional
analogue of the ``virtual surjection to pairs'' theorem
mentioned above:

\begin{theorem}[Virtual Surjections Theorem, Cohen--Shusterman {\cite[Theorem~1.3]{CohenShusterman2026}}]
  \label{thm:vsc}
  Let $m,n\in\N$ satisfy $2\leq m\leq n$, let
  $\Gamma_1,\Gamma_2,\dots,\Gamma_n$ be groups
  of type $\F_m$ (resp.\ $\FP_m$),
  let $\Gamma=\Gamma_1\times\Gamma_2\times\dots\times\Gamma_n$ be their direct product,
  let $G\leq \Gamma$
be a subgroup, and
assume that for all $i_1,i_2,\dots, i_m\in\{1,2,\dots,n\}$ 
with $i_1< i_2<\dots< i_m$ it holds that
the image of $G$ under the natural projection
$\Gamma\to\Gamma_{i_1}\times\dots\times\Gamma_{i_m}$ has finite index.
Then $G$ is of type $\F_m$ (resp.\ $\FP_m$).
\end{theorem}

Our proof of \cref{thm:main-fp} does not rely on any of the results
from \cite{CohenShusterman2026}, so
\cref{thm:main-fp} in conjunction with \cite[Theorem~3.10]{Kuckuck2014}
and \cite[Theorem~F]{KochloukovaLima2018}
provides an alternative proof of
\cref{thm:vsc}.

We currently do not know a counterexample to the statement in \cref{thm:main-fp}
when ``$Q$ is virtually nilpotent'' is replaced with ``$Q$ is of type $\F_{\min\{m,n\}+1}$''
(i.e., the condition in the general $n$-$(n+1)$-$(n+2)$ theorem),
nor do we know how to extend the proof beyond the virtually nilpotent setting.
A natural next broader class to consider 
for the quotients is the class of virtually polycyclic
groups. In this case, and more generally in the case where the
quotient has coherent integral group ring, we prove 
the weaker conclusion that the fibre product
has finitely generated integral homology up to the relevant degree:
\begin{athm}{theorem}{thm:main-coherent}[Weak symmetric $n$-$(n+1)$-$(n+2)$ theorem for
  fibre products over virtually polycyclic quotients]
Let $m,n,k,l\in\N_0$,
let
\[
\ses{
    1\to N_1\to\G_1\xrightarrow{\pi_1}Q\to1
} {
    1\to N_2\to\G_2\xrightarrow{\pi_2}Q\to1
}
\]
be short exact sequences of groups, assume that
$N_1$ is of type $\FP_k$, that $N_2$ is of type $\FP_l$,
that $\G_1$ is of type $\FP_m$, that $\G_2$ is of type $\FP_n$, 
and that $Q$ is virtually polycyclic (or, more generally, that $\Z Q$ is coherent). 
Then the fibre product
\[
  P=\{(\gamma_1,\gamma_2)\in\G_1\times\G_2
       \colon \pi_1(\gamma_1)=\pi_2(\gamma_2)\}
\]
satisfies
\[
\forall j\in\{0,1,\dots,\min\{m,n,k+l+1\}\}\colon 
\text{$H_j(P;\Z)$ is finitely generated.}
\]
\end{athm}

\paragraph{Proof idea.}
We will give a quick overview of the ideas behind the proof
of \cref{thm:main-fp}.
The proof proceeds by induction on the nilpotency class of $Q$, so
the main step concerns a single
central extension. Given short exact sequences
\[
\ses{
    1\to N_1\to\G_1\xrightarrow{\pi_1}Q\to1
} {
    1\to N_2\to\G_2\xrightarrow{\pi_2}Q\to1
}
\]
and a central extension
\begin{equation*}
  1\to Z\to Q\xrightarrow{p}\bar Q\to 1
\end{equation*}
one can consider the nested fibre products
\[\begin{aligned}
  \bar P&=\{(\gamma_1,\gamma_2)\in\G_1\times\G_2
       \colon (p\circ\pi_1)(\gamma_1)=(p\circ\pi_2)(\gamma_2)\}
       \\&\supseteq\{(\gamma_1,\gamma_2)\in\G_1\times\G_2
       \colon \pi_1(\gamma_1)=\pi_2(\gamma_2)\}=P.
\end{aligned}
\]
The first crucial observation is that $P$ is the kernel
of a homomorphism
\begin{equation*}
  \delta\colon\bar P\to Z,
  \ (\gamma_1,\gamma_2)
    \mapsto\pi_1(\gamma_1)\pi_2(\gamma_2)^{-1}.
\end{equation*}
That makes the problem of transferring finiteness properties from
$\bar P$ to $P$ amenable to methods from
the theory of Bieri--Neumann--Strebel--Renz $\Sigma$-invariants.
For every group one defines
the character sphere $S(G)=(\mathrm{Hom}(G,\R)\setminus\{0\})/{\sim}$
where two homomorphisms $\chi_1,\chi_2\colon G\to\R$
are $\sim$-equivalent if there is a $\lambda\in\R$ with
$\lambda>0$ such that $\chi_2=\lambda\chi_1$.
The Bieri--Neumann--Strebel--Renz invariants give for every group $G$
of type $\FP_m$
a chain of descending subsets
$S(G)\supseteq\Sig 1G\supseteq\Sig 2G\supseteq\dots$
such that
a subgroup $H\lhdnormal G$ with $G/H$ abelian
is of type $\FP_m$ if and only if
$[\chi]\in\Sig mG$ for all non-zero $\chi\colon G\to\R$ with
$\chi(H)=0$.
We prove the following lower bound for the
$\Sigma$-invariants of a fibre product (which
is known as Meinert's inequality in the case $Q=1$).
\begin{athm}{proposition}{prop:meinert-intro}
Let $k,l\in\N_0$, let
\[
\ses{
    1\to N_1\to\G_1\xrightarrow{q_1} Q\to1
} {
    1\to N_2\to\G_2\xrightarrow{q_2} Q\to1
}
\]
be short exact sequences of groups, let
\[
   P
  =\set{(\gamma_1,\gamma_2)\in\G_1\times\G_2\given
  q_1(\gamma_1)=q_2(\gamma_2)}
\]
be their fibre product, and assume that $ P$ is of type $\FP_{k+l+1}$.
Then $ N_1\times N_2\subseteq  P$ and
\begin{equation*}
  \set{[\chi]\in S(P)\given\text{
    $\chi|_{ N_1}\neq 0\neq \chi|_{ N_2}$,
  $[\chi|_{ N_1}]\in\Sigma^k( N_1;\Z)$,
  $[\chi|_{ N_2}]\in\Sigma^l( N_2;\Z)$}}
  \subseteq \Sigma^{k+l+1}( P;\Z).
\end{equation*}
\end{athm}
The proof combines
known properties of $\Sigma$-invariants but it requires
a certain
change of perspective. With the usual definition,
the invariant $\Sigma^m(G;\Z)$ is empty whenever $G$ is not 
of type $\FP_m$, so it only contains useful
information in dimensions up to the finiteness length of $G$ (i.e., the largest
$m\in\N_0\cup\{\infty\}$ such that $G$ is of type $\FP_m$).
But there is a well-known
alternative condition in terms of the vanishing of Novikov 
homology up to a certain degree
which is equivalent to the usual definition
in dimensions up to the finiteness length of $G$
but is not necessarily trivial beyond that range.
A substantial number of classical results from $\Sigma$-theory appear
to carry over to these higher-dimensional invariants
when one simply removes finiteness assumptions on the involved groups.
This allows one, e.g., to use such results in constructions involving
intermediate groups without making finiteness assumptions on these groups
(for example, in \cref{prop:meinert-intro} above, $N_1\times N_2$
need not satisfy any additional finiteness assumption beyond
what is implied by the assumptions on the factors).

We are certainly not the first to have noticed this
phenomenon. E.g., Papadima and Suciu~\cite[Section~8.2]{papadima2010bieri}
use the Novikov homology vanishing condition as a definition
and explicitly call the resulting invariants ``generalized $\Sigma$-invariants''
and Kochloukova and Mendonça~\cite{KochloukovaMendonca2022} use exactly the
kind of reasoning described above to get rid of finiteness assumptions
on intermediate groups (see, e.g., \cite[Proof of Theorem~1.5]{KochloukovaMendonca2022}).
The general principle still seems underexploited in the existing literature.

\paragraph{Structure of the article.} The rest of this article is organized as follows. 
In \cref{sec:fundamentals} we recall the required definitions 
and facts on finiteness
properties, Novikov homology, and $\Sigma$-invariants.
In \cref{sec:meinert} we prove the
Meinert inequality for fibre products, \cref{prop:meinert-intro} above as
\cref{prop:fibre-product-character}.
In \cref{sec:central-ext} we prove the key theorem on
fibre products over a central extension
and \cref{sec:nilpotent} completes the induction that
gives the main result. \Cref{thm:main-fp} is established
in \cref{cor:virtually-nilpotent-fp-2,cor:homotopical-virtually-nilpotent-fp}.
In \cref{sec:consequences}, we record the
Virtual Surjections Theorem as a consequence and prove that the 
term $k+l+1$ in the conclusion of \cref{thm:main-fp} is optimal.
Finally, in \cref{sec:coherent-quotients}, we prove 
the weak result for quotients with coherent integral group ring, \cref{thm:main-coherent},
as \cref{thm:coherent-quotient-homology}.

\section{Preliminaries}
\label{sec:fundamentals}

\subsection{Finiteness properties and \texorpdfstring{$\Sigma$}{Sigma}-invariants}

We recall the finiteness properties and the homological
Bieri--Neumann--Strebel--Renz invariants used in the sequel.

\begin{athm}{definition}{def:finiteness-properties}[Finiteness properties]
Let $m\in \N_0$ and let $R$ be a ring. A left $R$-module $M$ is said to be of type $\FP_m$ if it admits a projective resolution
\[
  \cdots\to P_1\to P_0\to M\to0
\]
in which $P_i$ is finitely generated for all $i\in\{0,1,\dots,m\}$.
A left $R$-module $M$ is said to be of type $\FP_\infty$ if it is of type
$\FP_n$ for all $n\in\N_0$.
A group $G$ is said to be of type $\FP_m$ (resp.\ $\FP_\infty$) if $\Z$, regarded as
a trivial left $\Z G$-module, is of type $\FP_m$ (resp.\ $\FP_\infty$).

A group $G$ is said to be of type $\F_m$ if it admits an Eilenberg--Mac Lane complex
(i.e., a CW complex with $G$ as fundamental group and contractible universal cover) whose $m$-skeleton is finite. A group is said to be of type $\F_\infty$ if it is of type
$\F_n$ for all $n\in\N_0$.
\end{athm}

We will only need a few standard finiteness facts:

\begin{athm}{lemma}{lem:finiteness-facts}
Let $m\in\N_0\cup\{\infty\}$. Then the following hold:
\begin{enumerate}[(i)]
\item \label{it:fp1} For every group $G$ it holds that
$G$ is of type $\F_1$ if and only if $G$ is of type $\FP_1$ if and only if $G$ is finitely generated.
\item \label{it:wall} A group is of type $\F_{m+2}$ if and only if it is finitely presented
and of type $\FP_{m+2}$.
\item \label{it:finite-index} Type $\FP_m$ and type $\F_m$ are invariant under passage to finite-index
subgroups and finite-index overgroups.
\item \label{it:ext-preserve-fp} For every short exact sequence $1\to A\to B\to C\to1$ of groups with 
$A$ of type $\FP_m$ (resp.\ $\F_m$) and $C$ of
type $\FP_m$ (resp.\ $\F_m$) it holds that $B$ is of type $\FP_m$ (resp.\ $\F_m$).
\item \label{it:nilpotent-fpinfty} If $Q$ is a finitely generated nilpotent group, then $Q$ and every subgroup of
$Q$ are of type $\F_\infty$.
\end{enumerate}
\end{athm}
\begin{aproof}
  \Cref{it:fp1} is elementary.
\Cref{it:finite-index} can be found, for the $\FP_m$ properties, e.g., in
Brown~\cite[Chapter~VIII, Proposition~5.1]{Brown1982} and, for the $\F_m$ properties,
e.g., in Geoghegan~\cite[Corollary~7.2.4]{Geoghegan2008}.
\Cref{it:wall} is Wall's criterion~\cite{Wall1965,Wall1966} (see also
Brown~\cite[Section~VIII.7]{Brown1982}).
\Cref{it:ext-preserve-fp} can be proved, for the $\FP_m$ properties,
using the Bieri--Eckmann criterion (see Bieri \& Eckmann~\cite{BieriEckmann1974})
and a Lyndon--Hochschild--Serre spectral sequence
(see, e.g., Bieri~\cite[Proposition~2.7 and the following exercise]{Bieri1976}) and, 
for the $\F_m$ properties,
using Wall's criterion or direct construction (see, e.g., Geoghegan~\cite[Chapter~7]{Geoghegan2008}).
\Cref{it:nilpotent-fpinfty} is a consequence of \cref{it:ext-preserve-fp} and the fact that finitely generated nilpotent groups and all their subgroups are polycyclic.
\end{aproof}

We next recall the definition
of the $\Sigma$-invariants due to Bieri, Neumann, Strebel, and Renz~\cite{BieriNeumannStrebel1987,BieriRenz1988}.

\begin{athm}{definition}{def:sigma-invariants}[\texorpdfstring{$\Sigma$}{Sigma}-invariants]
Let $G$ be a group. For every character $\chi\colon G\to\R$ we denote
$[\chi]=\set{\lambda \chi\given \lambda\in\R,\,\lambda>0}$.
Moreover, we write
\[
  S(G)=\set{[\chi]\given \chi\in\mathrm{Hom}(G,\R),\,\chi\neq 0}
\]
and we call $S(G)$ the character sphere of $G$.
For every subgroup $H\leq G$ we denote
\begin{equation*}
  S(G,H) = \set{[\chi]\in S(G)\given \chi(H)=0}.
\end{equation*}
In addition, for every character $\chi\colon G\to\R$, we denote
\[
  G_\chi=\{g\in G\colon \chi(g)\ge0\}.
\]
For $m\in\N_0$ and $A$ a left $\Z G$-module
the homological invariant $\Sigma^m(G;A)$ is defined to be the set of
classes $[\chi]\in S(G)$ such that $A$ is of type $\FP_m$ as a left
$\Z G_\chi$-module.\footnote{In particular, note that $\Sigma^0(G;\Z)=S(G)$ for every group $G$.}
\end{athm}

\begin{remark}
  Note that we allow the notation $[\chi]$ also when $\chi$ is the zero character,
  but $[\chi]\in S(G)$ implies that $\chi\neq 0$.%
\end{remark}
\begin{remark}
  We will only consider the $\Sigma$-invariants
  $\Sigma^m(G;\Z)$ with integral coefficients here.
  We retain the coefficient module in the notation
  because $\Sigma^m(G)$ usually denotes
  a different, homotopical version of the invariant.
\end{remark}

The $\Sigma$-invariants of a group $G$ 
determine the finiteness properties $\FP_m$ of kernels of maps from $G$ 
to abelian groups.

\begin{athm}{theorem}{thm:renz-fp}[Bieri--Renz criterion, {\cite[Theorem~B]{BieriRenz1988}}]
Let $m\in\N_0$.  Let $G$ be a group of type $\FP_m$, and let $N\lhdnormal G$ be a
normal subgroup such that $G/N$ is abelian.  Then $N$ is of type $\FP_m$ if and
only if
\[
  S(G,N)=\{[\chi]\in S(G)\colon \chi(N)=0\}
  \subseteq \Sigma^m(G;\Z).
\]
\end{athm}

\subsection{Novikov--Sikorav completions}

For a group $G$, we will in the following write functions $a=(a_g)_{g\in G}\colon G\to \Z$ as formal
infinite series $\sum_{g\in G} a_g g$. The set $\Z^G$ of such functions naturally forms an abelian group
under pointwise addition and can further be endowed with a partial multiplication
(extending the multiplication of the group ring $\Z G$), defined for all $a,b\in\Z^G$ such that
it holds for all $x\in G$ that $\#\{g\in G\colon a_g b_{g^{-1}x}\neq 0\}<\infty$ by
\begin{equation}
  \textstyle
    ab
    =
    \pr*{\shortsum_{g\in G}a_g g}
    \pr*{\shortsum_{h\in G}b_h h}
    =
    \shortsum_{x\in G}\pr*{\shortsum_{g\in G}a_g b_{g^{-1}x}}x.
\end{equation}
This multiplication is not necessarily associative.

The homological $\Sigma$-invariants can be characterized in terms of a
particular subring of $\Z^G$, introduced by Novikov~\cite{Novikov1981}
(in the abelian case) and Sikorav~\cite[Chapter~IV]{Sikorav1987}.

\begin{athm}{definition-and-lemma}{def:novikov}[Novikov--Sikorav completion]
Let $G$ be a group and let $\chi\colon G\to\R$ be a character.  The
Novikov--Sikorav completion of $\Z G$ in the direction $\chi$ is
\begin{equation}\label{eq:novdef}
  \textstyle
  \NS{G}{\chi}
  =
  \set*{
  \shortsum_{g\in G} a_g g \in \Z^G\given
  [\forall C\in\R\colon
  \#\{g\in G\colon a_g\ne 0,\, \chi(g)\le C\}<\infty]
  }.
\end{equation}
For all $a,b\in\NS G\chi$ the product $ab$ is defined
and lies in $\NS G\chi$
and this multiplication is associative on $\NS G\chi$. 
Therefore, $\NS G\chi$ is an associative unital ring.
\end{athm}

The group ring $\Z G$ is a subring of $\NS{G}{\chi}$ and $\NS G\chi$ thus
naturally has the structure of a $\Z G$-(bi)module. Considering it as a right $\Z G$-module, we
can therefore consider for every $i\in\N_0$ the
homology groups
\[
  H_i(G;\NS G\chi)
  =\operatorname{Tor}_i^{\Z G}(\NS G\chi,\Z),
\]
where $\Z$ is the trivial left $\Z G$-module.
These are closely tied to the homological $\Sigma$-invariants via the following criterion, 
which is 
due to Sikorav~\cite[Chapter~IV, Section~2.5]{Sikorav1987} in the $1$-dimensional case
and established by 
Schweitzer in the general high-dimensional form (see Bieri~\cite[Theorem~A.1]{Bieri2007}).
\begin{athm}{theorem}{thm:novikov}[Novikov--Sikorav criterion]
Let $m\ge1$, let $G$ be a group of type $\FP_m$, and let $\chi\colon G\to\R$ be a
non-zero character.  Then the following are equivalent:
\begin{enumerate}[(i)]
  \item \label{it:sigma} It holds that $[\chi]\in\Sigma^m(G;\Z)$.
  \item \label{it:novvanish} For all $i\in\{0,1,\dots,m\}$ it holds that $H_i(G;\NS{G}{\chi})=0$.
\end{enumerate}
\end{athm}

\begin{remark}
  Sign conventions in the literature unfortunately differ. Sometimes the 
  inequality in \cref{eq:novdef} is reversed.
  In that case, the Novikov--Sikorav criterion has
  $-\chi$ instead of $\chi$ in either \cref{it:novvanish} or \cref{it:sigma}.
\end{remark}

\begin{athm}{lemma}{lem:invert}
Let $G$ be a group, let $\chi\colon G\to\R$ be a homomorphism, and let $g\in G$ with
$\chi(g)\ne0$.  Then $g-1$ is invertible in
$\NS{G}{\chi}$.
\end{athm}
\begin{aproof}
  Note that in $\Z^G$ it holds that 
  \begin{equation*}
    \textstyle(g-1)\pr[\big]{\sum_{j\ge0}-g^j} = 1 = \pr[\big]{\sum_{j\ge0}-g^j}(g-1)
    \;\;\text{and}\;\;
    (g-1)\pr[\big]{\sum_{j\ge1}g^{-j}} = 1 = \pr[\big]{\sum_{j\ge1}g^{-j}}(g-1).
  \end{equation*}
If $\chi(g)>0$, then $\sum_{j\ge0}-g^j\in\NS G\chi$ and if
$\chi(g)<0$ then $\sum_{j\ge1}g^{-j}\in\NS G\chi$.
\end{aproof}

\begin{athm}{remark}{rem:novikov-finiteness-assumption}
  It is worth noting that with our definition of the $\Sigma$-invariants, 
  if for some $m\in\N_0$ a group $G$ has $[\chi]\in\Sigma^m(G;\Z)$
  for a character $\chi\colon G\to\R$,
  then $G$ is automatically\footnote{This uses that for any
  $t\in G$ with $\chi(t)<0$
  it holds that $\Z G$ is an ascending union $\bigcup_{n\in\N} t^n\Z G_\chi$,
  hence flat as a right $\Z G_\chi$-module and $\Z G\tensor_{\Z G_\chi}\Z=\Z$;
  so a resolution of $\Z$ as 
  a $\Z G_\chi$ module through extension of scalars yields a resolution
  of $\Z$ as a $\Z G$-module and this preserves 
  projectivity and finite generation.} of type $\FP_m$. In fact, this assumption is often
  added redundantly to the definition of the $\Sigma$-invariants.
  
  By contrast, condition (ii) in \cref{thm:novikov} above does not usually imply
  that $G$ is of type $\FP_m$. 
  For an example, let $K$ be any group,
  take $G=K\times\Z=K\times\ang{t}$, and let 
  $\chi\colon G\to\Z$ be projection onto the second factor. 
  Now $t$ is central and $t-1$ is invertible in 
  $\NS G\chi$ by \cref{lem:invert}.
  Central elements act trivially on group homology, 
  so multiplication by $t-1$ is the
  zero map on $H_*(G;\NS G\chi)$;
  but, since $t-1$ is invertible, this map is also an automorphism. 
  Hence all Novikov homology of $G$ vanishes.
  But if $K$ is not finitely generated, neither is $G$.
  
  Therefore, the assumption in \cref{thm:novikov} that
  $G$ be of type $\FP_m$ is necessary: The equivalence does not usually hold beyond
  the finiteness range of $G$. This should not be seen as a deficiency. It means that
  Novikov acyclicity can still detect some form of the controlled directional
  part of the $\Sigma$-condition without imposing ambient finiteness.
\end{athm}

\begin{athm}{definition}{def:extended-sigma}[Novikov-$\Sigma$-invariants]
  Let $G$ be a group and $m\in\N_0$. We denote 
  \begin{equation*}
    \Nov{m}{G}
    =
    \set[\big]{[\chi]\in S(G)\given \br[\big]{\forall i\in\{0,1,\dots,m\}\colon H_i(G;\NS{G}{\chi})=0}}.
  \end{equation*}
\end{athm}

\begin{athm}{lemma}{lem:novikov-alternative-form}[Novikov--Sikorav criterion, alternative form]
Let $m\in\N_0$ and let $G$ be a group of type $\FP_m$.  Then
\begin{equation*}
  \Nov{m}{G}=\Sigma^m(G;\Z).
\end{equation*}
\end{athm}
\begin{aproof}
For $m\geq1$ this is \cref{thm:novikov}. For $m=0$,
note that \cref{lem:invert}
shows that for every $[\chi]\in S(G)$, $x\in\NS G\chi$, $g\in G$ with $\chi(g)\neq 0$
it holds that $x=x(g-1)^{-1}(g-1)$, hence $x\in (\NS G\chi)(G-1)$. Therefore,
$H_0(G;\NS G\chi)=(\NS G\chi)_G=0$ and so $\Nov{0}{G}=S(G)=\Sigma^0(G;\Z)$.
\end{aproof}

\begin{remark}
  Note that, with our conventions, the inclusion
  $\Sigma^m(G;\Z)\subseteq \Nov{m}{G}$ is in fact
  satisfied for all groups $G$ and all $m\in\N_0$, since
  $\Sigma^m(G;\Z)=\emptyset$
  when $G$ is not of type $\FP_m$.
\end{remark}

We will frequently use a characterization of $\Nov{m}{G}$ that a priori looks
stronger than the defining condition, see \cref{lem:restriction-vanishing}
below. This is essentially due to Schütz~\cite[Corollary~2.5]{Schuetz2009}. It can be proved
by a straightforward application of the universal coefficient spectral sequence, but we 
include a direct proof here for the reader's convenience.

First, we will frequently need the following observation,
which essentially follows immediately from the definitions.
\begin{athm}{lemma}{lem:restriction-module}
Let $H\le G$, let $\chi\colon G\to\R$ be a character, let
$\psi=\chi|_H$, and assume $\psi\neq 0$.  Then the inclusion 
$\Z^H\to\Z^G,\,\sum_{h\in H}a_h h\mapsto \sum_{h\in H}a_h h$
restricts to an inclusion $\NS{H}{\psi}\to\NS{G}{\chi}$ making
$\NS H\psi$ a subring of $\NS G\chi$. Consequently, every 
$\NS G\chi$-module can be naturally regarded as a $\NS H\psi$-module with the 
$\NS H\psi$-module structure extending the $\Z H$-module structure that is conferred
via the inclusion $\Z H\to \Z G\to \NS G\chi$.
\end{athm}

\begin{athm}{lemma}{lem:restriction-vanishing-general}
Let $G$ be a group, let $\mathbb ZG\to S$ be a ring homomorphism, 
let $M$ be a right $S$-module and let $m\in\N_0$
satisfy for all $i\in\N_0$ with $i\leq m$ that
$H_i(G;S)=0$.
Then it holds for all $i\in\N_0$ with $i\leq m$ that
\[H_i(G;M)=0.\]
\end{athm}
\begin{aproof}
Let $P_*\to\Z$ be a projective resolution of the trivial left
$\Z G$-module $\Z$, let $(C_*,d_*)=S\otimes_{\Z G}P_*$, and
for all $i\in\N_0$ denote $B_i=\im d_{i+1}$.
Note that for all $i\in \N_0$ it holds that
\begin{equation}\label{eq:4}
  M\otimes_{S}C_i
  =
  M\otimes_{S}(S\otimes_{\Z G}P_i)
  \cong
  M\otimes_{\Z G}P_i.
\end{equation}
We thus have for all $i\in\N_0$ that
\begin{equation*}
  H_i(G;M)=H_i(M\otimes_{\Z G}P_*)=H_i(M\otimes_S C_*).
\end{equation*}
For all $i\in\N_0$ with $i\leq m$ the assumption that $H_i(G;S)=0$ implies that
\begin{equation}\label{eq:33}
  H_i(C_*)=H_i(G;S)=0.
\end{equation}
Since each $P_i$ is a direct summand of a free $\Z G$-module, extension of
scalars makes $C_i$ a direct summand of a free $S$-module.  Thus
every $C_i$ is projective over $S$.
By \cref{eq:33}, we have for every $i\in\N$ with $i\leq m$ a short exact sequence
\begin{equation} \label{eq:ses}
  0\to B_i\to C_i\xrightarrow{d_i}B_{i-1}\to0.
\end{equation}
Since $H_0(C_*)=0$, we have that $B_0=C_0$ is projective and hence, if $m\geq 1$, the
sequence \cref{eq:ses} splits for $i=1$.
Next, note that for all $i\in\N$ with $i\leq m-1$ where \cref{eq:ses} splits,
it holds that $B_i$ is projective as a direct summand of the projective
module $C_i$, hence
the sequence
\begin{equation*} 
  0\to B_{i+1}\to C_{i+1}\xrightarrow{d_{i+1}}B_{i}\to0
\end{equation*}
likewise splits.
Induction therefore establishes 
that for all $i\in\N$ with $i\leq m$ the sequence \cref{eq:ses} is split.
For every $i\in\N$ with $i\leq m$ we therefore have the following diagram with an 
exact middle row where the
vertical arrow on the left is a surjection by right-exactness of tensoring and the 
vertical arrow on the right is an injection because the inclusion map
$B_{i-1}\to C_{i-1}$ is the identity for $i=1$ and a split injection for $2\leq i\leq m$:
\[
\begin{tikzcd}[column sep=large, row sep=large]
  M\otimes_S C_{i+1}
    \arrow[d, two heads]
    \arrow[dr, "\id\otimes d_{i+1}"]
  \\
  M\otimes_S B_i
    \arrow[r, hook]
  &
  M\otimes_S C_i
    \arrow[r, two heads]
    \arrow[dr, "\id\otimes d_{i}"]
  &
  M\otimes_S B_{i-1}
    \arrow[d, hook]
  \\
  &&M\otimes_S C_{i-1}
\end{tikzcd}
\]
This shows for all 
$i\in\N$ with $i\leq m$ 
that $\ker(\id_M\tensor d_i)=M\tensor_S B_i=\im(\id\tensor d_{i+1})$
and hence
\begin{equation}\label{eq:55}
  H_i(M\tensor_{S}C_*)=0.
\end{equation}
In addition, $d_1\colon C_1\to B_0=C_0$ is surjective, so
$\id\tensor d_1\colon M\otimes_S C_1\twoheadrightarrow M\otimes_S C_0$
is as well, and therefore $H_0(M\otimes_S C_*)=0$.
Combining this and \cref{eq:55} with \cref{eq:4} finishes the proof.
\end{aproof}

\begin{athm}{corollary}{lem:restriction-vanishing}
Let $G$ be a group, let $\chi\colon G\to\R$ be a non-zero character,
and let $m\in\N_0$. Then the following are equivalent:
\begin{enumerate}[(i)]
  \item \llabel{it:1} It holds that $[\chi]\in\Nov{m}{G}$.
  \item \llabel{it:2} For every right $\NS G\chi$-module $M$ 
  and every $i\in\{0,1,\dots,m\}$ it holds that $H_i(G;M)=0$.
\end{enumerate}
\end{athm}
\begin{aproof}
  The implication \lref*{it:2} $\implies$ \lref*{it:1} is simply the case
  $M=\NS G\chi$. The implication 
  \lref*{it:1} $\implies$ \lref*{it:2} is \cref{lem:restriction-vanishing-general}
  applied with $S=\NS G\chi$.
\end{aproof}

\section{Meinert's inequality for fibre products}
\label{sec:meinert}

The following criterion is known as Meinert's inequality
when the Novikov-$\Sigma$-invariants
are replaced by regular $\Sigma$-invariants and
$A\times B$ is additionally assumed to be of type $\FP_{k+l+1}$
(unpublished by Meinert, but recorded in Gehrke~\cite[Lemma~9.1]{Gehrke1998}).
Our version is essentially due to Schütz~\cite[Corollary~3.3]{schuetz2008} but (crucially for the later
application) does not make any finiteness assumptions on the groups involved.
We include a proof here for the reader's convenience.

\begin{athm}{lemma}{lem:direct-product-add-up}
Let $k,l\in\N_0$, let $A$ and $B$ be groups, let $G=A\times B$, and let
$\chi\colon G\to\R$ satisfy
\[
  [\chi|_A]\in\Nov{k}{A}
  \qquad\text{and}\qquad
  [\chi|_B]\in\Nov{l}{B}.
\]
Then
\[
  [\chi]\in\Nov{k+l+1}{G}.
\]
\end{athm}
\begin{aproof}
Throughout this proof, denote $M=\NS G\chi$,
$\psi=\chi|_A$, and $\eta=\chi|_B$.
Note that the assumption that $[\chi|_A]\in\Nov kA$ implies that $\chi|_A\neq 0$ and hence $\chi\neq 0$.
Consider the Lyndon--Hochschild--Serre spectral sequence for
the canonical short exact sequence
\[
  1\to A\to A\times B\to B\to1
\]
with coefficients in the right $\Z G$-module $M$:
\begin{equation}
  \llabel{eq:lhs}
  E^2_{p,q}=H_p(B;H_q(A;M))\Longrightarrow H_{p+q}(G;M).
\end{equation}
Observe that \cref{lem:restriction-module}
shows that $M$ has a right $\NS A\psi$-module structure extending the right $\Z A$-module structure conferred by the inclusion $\Z A\to \Z G\to \NS G\chi$.
The assumption that $[\psi]\in\Nov{k}{A}$ and \cref{lem:restriction-vanishing} hence show that
\begin{equation}\llabel{eq:HqAM}
  \forall q\in\{0,1,\dots, k\}\colon H_q(A;M)=0.
\end{equation}

Recall the definition of the quotient action in the Lyndon--Hochschild--Serre
spectral sequence \lref{eq:lhs}: 
Choose a free left
$\Z A$-resolution $P_*\xrightarrow{\varepsilon}\Z$ and for every $b\in B$ a lift $\tilde b\in G$.
Then choose for every $b\in B$ a chain map $\phi^b_*\colon P_*\to P_*$ satisfying for all
$i\in\N_0$, $\gamma\in P_i$, $a\in A$ that
$\phi^b_i(a\gamma)=\tilde b^{-1}a\tilde b\phi^b_i(\gamma)$
and $\varepsilon\circ\phi^b_0=\varepsilon$.
Then there is
for every $b\in B$ a chain map $\psi^b_*\colon M\tensor_{\Z A} P_*\to M\tensor_{\Z A}P_*$
which for all $m\in M$, $i\in\N_0$, $\gamma\in P_i$ satisfies
$\psi^b_i(m\tensor \gamma)=(m\tilde b)\tensor\phi^b_i(\gamma)$. The right multiplication
by $b\in B$ on $H_*(A;M)$ is simply the map induced by $\psi^b_*$ on homology.
In the present situation every $b\in B$ has the canonical lift
$\tilde b=(1,b)\in G$, which centralizes $A$, so the
map $\phi^b_*$ above can be chosen to be the identity chain map.
Thus, in \lref{eq:lhs}, the action of $B$ on $H_q(A;M)$
is induced by the unique action on $M\tensor_{\Z A}P_*$
which satisfies for all $i\in\N_0$, $m\in M$, $\gamma\in P_i$, $b\in B$ that
\begin{equation}\llabel{eq:C-action-on-A-chains}
  (m\tensor \gamma)b=(mb)\tensor \gamma.
\end{equation}
We now show that this action can be extended to an action of $\NS B\eta$ on $H_*(A;M)$.
By \cref{lem:restriction-module}, $M$
is naturally a right $\NS{B}{\eta}$-module. Since $\NS B\eta$ centralizes $\Z A$,
the action 
of $\NS B\eta$ commutes with the action of $\Z A$.
Therefore, it holds for every $\lambda\in\NS B\eta$ that the right multiplication map
$M\to M,\,m\mapsto m\lambda$ is a homomorphism of right $\Z A$-modules.
Functoriality of $\Tor_*^{\mathbb ZA}(-,\mathbb Z)$
therefore equips
\[
H_*(A;M)=\Tor_*^{\mathbb ZA}(M,\mathbb Z)
\]
with a natural right $\NS B\eta$-module structure. 
Explicitly, if $P_*\to\mathbb Z$ is any free left $\mathbb ZA$-resolution, this action 
on $H_*(A;M)$ is induced by the unique action of $\NS B\eta$ on $M\tensor_{\Z A}P_*$ which
satisfies for all $i\in\N_0$, $m\in M$, $\gamma\in P_i$, $\lambda\in \NS B\eta$ that
\[
(m\otimes \gamma)\lambda=(m\lambda)\otimes \gamma,
\]
and so it extends the $B$-action on $H_*(A;M)$ in \lref{eq:C-action-on-A-chains}.
The assumption that
$[\eta]\in\Nov{l}{B}$ and \cref{lem:restriction-vanishing} now give
\[
  \forall p\in\{0,1,\dots,l\}\colon H_p(B;H_q(A;M))=0.
\]
Combining this with \lref{eq:HqAM} shows that $E^2_{p,q}=0$ whenever
$p\leq l$ or $q\leq k$.  If $p+q\leq k+l+1$, at least one of these
inequalities holds.  Therefore
\begin{equation}
  \forall p,q\in\N_0\colon p+q\leq k+l+1\implies E^2_{p,q}=0.
\end{equation}
This proves
\begin{equation*}
  \forall s\in\{0,1,\dots,k+l+1\}\colon H_s(G;M)=0
\end{equation*}
and therefore $[\chi]\in\Nov{k+l+1}{G}$.
\end{aproof}

The following normal-subgroup ascent property was observed
before by Kochloukova and
Mendon{\c c}a (cf.\ \cite[Proof of Theorem~1.5]{KochloukovaMendonca2022}).

\begin{athm}{lemma}{lem:normal-subgroup-novikov-ascent}
Let $m\in\N_0$, let $G$ be a group, let $N\lhdnormal G$ be a normal subgroup, 
and let
$\chi\colon G\to\R$ satisfy $[\chi|_N]\in\Nov{m}{N}$. Then
$
  [\chi]\in\Nov{m}{G}
$.
\end{athm}
\begin{aproof}
  Let $\psi=\chi|_N$. By
\cref{lem:restriction-module}, restriction of scalars makes $\NS G\chi$ a right
$\NS N\psi$-module.  The assumption that $[\psi]\in\Nov{m}{N}$ and
\cref{lem:restriction-vanishing} therefore give
\[
\forall q\in\{0,1,\dots,m\}\colon 
  H_q(N;\NS G\chi)=0.
\]
The Lyndon--Hochschild--Serre spectral sequence for
\[
  1\to N\to G\to G/N\to1
\]
with coefficients in the right $\Z G$-module $\NS G\chi$ has the form
\[
  E^2_{p,q}=H_p(G/N;H_q(N;\NS G\chi))\Longrightarrow H_{p+q}(G;\NS G\chi).
\]
Its terms in total degree at most $m$ vanish, and hence
\[
\forall s\in\{0,1,\dots,m\}\colon
  H_s(G;\NS G\chi)=0\qedhere
\]
\end{aproof}

Simply combining the two lemmas above immediately gives the
following strengthening of the usual Meinert's inequality
where the ambient finiteness is not necessarily provided by the direct product itself
but by a larger overgroup.

\begingroup
\begin{athm}{proposition}{prop:meinert-split}
  Let $k,l\in\N_0$, let $G$ be a group, and let
  $A,B\leq G$ be subgroups
  such that $[A,B]=1$, $A\cap B=1$ and $AB\cong A\times B$ is normal
  (or subnormal) in $G$. Then
  \begin{equation*}
    \set{[\chi]\in S(G)\colon \text{$[\chi|_A]\in\Nov kA$ and $[\chi|_B]\in\Nov lB$}}\subseteq\Nov{k+l+1}G.
  \end{equation*}
  If, in addition, $G$ is of type $\FP_{k+l+1}$, the
  same statement holds with $\Sigma_{\mathrm{Nov}}^\cdot(-;\Z)$ replaced by $\Sigma^\cdot(-;\Z)$ throughout.
\end{athm}
\begin{aproof}
  For any 
  $[\chi]\in S(G)$ with $[\chi|_A]\in\Nov kA$ and $[\chi|_B]\in\Nov lB$
  \cref{lem:direct-product-add-up} implies that
  $[\chi|_{AB}]\in\Nov{k+l+1}{AB}$
and then applying \cref{lem:normal-subgroup-novikov-ascent} (once
if $AB$ is normal in $G$ or several times if $AB$ is subnormal)
yields $[\chi]\in\Nov{k+l+1}{G}$. 
The analogous statement for the regular $\Sigma$-invariants
follows by the Novikov--Sikorav criterion,
\cref{lem:novikov-alternative-form}, if $G$ is of type
$\FP_{k+l+1}$ (see also the remark following \cref{lem:novikov-alternative-form}
for why no finiteness assumptions on $A$ and $B$ are necessary).
\end{aproof}
\endgroup

\begin{athm}{proposition}{prop:fibre-product-character}
Let $k,l\in\N_0$, let
\[
\ses{
    1\to\bar N_1\to\G_1\xrightarrow{q_1}\bar Q\to1
} {
    1\to\bar N_2\to\G_2\xrightarrow{q_2}\bar Q\to1
}
\]
be short exact sequences of groups, let
\[
  \bar P
  =\set{(\gamma_1,\gamma_2)\in\G_1\times\G_2\given
  q_1(\gamma_1)=q_2(\gamma_2)}
\]
be their fibre product, and let
$\chi\colon \bar P\to\R$ satisfy $[\chi|_{\bar N_1}]\in\Nov{k}{\bar N_1}$
and $[\chi|_{\bar N_2}]\in\Nov{l}{\bar N_2}$ (identifying
$\bar N_1$ and $\bar N_2$ with the coordinate
subgroups $\bar N_1\times1$ and $1\times\bar N_2$ of $\bar P$).
Then $[\chi]\in \Nov{k+l+1}{\bar P}$.

  If, in addition, $\bar P$ is of type $\FP_{k+l+1}$, the
  same statement holds with $\Sigma_{\mathrm{Nov}}^\cdot(-;\Z)$ replaced by $\Sigma^\cdot(-;\Z)$ throughout.
\end{athm}
\begin{aproof}
  This immediately follows from \cref{prop:meinert-split}
  after observing that $\bar N_1\times \bar N_2$
  is the kernel of the homomorphism
$\rho\colon\bar P\to\bar Q$, $(\gamma_1,\gamma_2)\mapsto q_1(\gamma_1)=q_2(\gamma_2)$.
\end{aproof}

\section{The central-extension step}
\label{sec:central-ext}

We now prove the induction step over central extensions.

\begin{athm}{lemma}{lem:abelian-quotient-fibre-product}
Let $k,l\in\N_0$, let
\[
\ses{
    1\to\bar N_1\to\G_1\xrightarrow{q_1}\bar Q\to1
} {
    1\to\bar N_2\to\G_2\xrightarrow{q_2}\bar Q\to1
}
\]
be short exact sequences of groups, let
\[
  \bar P
  =\set{(\gamma_1,\gamma_2)\in\G_1\times\G_2\given
  q_1(\gamma_1)=q_2(\gamma_2)}
\]
be their fibre product,
let $A$ be a finitely generated abelian group, let
$\delta\colon\bar P\to A$ be a homomorphism,
assume that $\delta(\bar N_1)=A=\delta(\bar N_2)$,
for $i\in\{1,2\}$ let
\[
  P=\ker\delta
  \qquad\text{and}\qquad
  N_i=\ker(\delta|_{\bar N_i}),
\]
and assume that $N_1$ is of type $\FP_k$ and $N_2$ is of type $\FP_l$.
Then
\begin{enumerate}[label=(\roman*)]
  \item \llabel{it:factor-finiteness} the groups $\bar N_1$ and $\bar N_2$
  are of types $\FP_k$ and $\FP_l$, respectively,
  \item \llabel{it:novikov}
  $S(\bar P,P)\subseteq\Nov{k+l+1}{\bar P}$, and
  \item \llabel{it:finiteness} for every $d\in\N_0$ such that
  $\bar P$ is of type $\FP_d$ it holds that
  $P$ is of type $\FP_{\min\{d,k+l+1\}}$.
\end{enumerate}
\end{athm}

\begin{aproof}
  Throughout this proof, for $i\in\{1,2\}$
  let $\delta_i=\delta|_{\bar N_i}$.
For $i\in\{1,2\}$, there is a short exact sequence
\[
  1\to N_i\to\bar N_i\xrightarrow{\delta_i}\im \delta_i\to1.
\]
The assumption that $N_1$ is of type $\FP_k$
and $N_2$ is of type $\FP_l$,
the assumption that $A$ is finitely generated,
and \cref{it:ext-preserve-fp,it:nilpotent-fpinfty}
in \cref{lem:finiteness-facts}
therefore
prove \lref{it:factor-finiteness}.
Combining \lref{it:factor-finiteness}
with the Bieri--Renz criterion, \cref{thm:renz-fp}, establishes that
\begin{equation}
  \llabel{eq:factors}
  S(\bar N_1,N_1)\subseteq \Sigma^k(\bar N_1;\Z)
  \qquad\text{and}\qquad
  S(\bar N_2,N_2)\subseteq \Sigma^l(\bar N_2;\Z).
\end{equation}
Note that the assumption that $\delta(\bar N_1)=A$
implies that $\delta$ is an epimorphism. 
This and the assumption that $\ker\delta = P$
ensure that 
for every non-zero character $\chi\colon\bar P\to \R$ with
$\chi(P)=0$ there is a non-zero character $\lambda\colon A\to\R$ such 
that $\chi=\lambda\circ\delta$. 
The assumption that $\delta(\bar N_1)=A=\delta(\bar N_2)$
therefore ensures that for every
non-zero character $\chi\colon\bar P\to \R$ with
$\chi(P)=0$ the restrictions
$\chi|_{\bar N_1}$ and $\chi|_{\bar N_2}$ are
also non-zero.
Combining this with
\lref{eq:factors} shows
for every non-zero character $\chi\colon \bar P\to\R$ with $\chi(P)=0$ and hence
$\chi(N_1)=0=\chi(N_2)$
that
\begin{equation}
  [\chi|_{\bar N_1}]\in \Sigma^k(\bar N_1;\Z)=\Nov k{\bar N_1}
  \qquad\text{and}\qquad
  [\chi|_{\bar N_2}]\in \Sigma^l(\bar N_2;\Z)=\Nov l{\bar N_2}.
\end{equation}
\cref{prop:fibre-product-character} hence
establishes \lref{it:novikov}.

\Lref{it:novikov} and the Novikov--Sikorav criterion, \cref{lem:novikov-alternative-form}, 
now show that for all $d\in\N_0$ with $d\leq k+l+1$ such that $\bar P$ is of type
$\FP_d$ it holds that
\[
  S(\bar P,P)
  \subseteq\Nov{k+l+1}{\bar P}
  \subseteq\Nov{d}{\bar P}
  =\Sigma^d(\bar P;\Z).
\]
The Bieri--Renz criterion, \cref{thm:renz-fp}, hence establishes \lref{it:finiteness}.
\end{aproof}

\begin{athm}{proposition}{prop:central-step}
Let $k,l\in\N_0$, let
\[
\ses{
  	  1\to N_1\to \G_1\xrightarrow{\pi_1}Q\to1
} {
    1\to N_2\to \G_2\xrightarrow{\pi_2}Q\to1
}
\]
be short exact sequences of groups, let
\begin{equation}\llabel{eq:centr}
  1\to Z\to Q\xrightarrow{p}\bar Q\to1
\end{equation}
be a central extension of groups with $Z$ finitely generated,
let
\begin{equation*}\begin{aligned}
  \bar P&=
  \{(\gamma_1,\gamma_2)\in\G_1\times\G_2\colon
  p\pi_1(\gamma_1)=p\pi_2(\gamma_2)\},
  \\
  P&=
  \{(\gamma_1,\gamma_2)\in\G_1\times\G_2\colon
  \pi_1(\gamma_1)=\pi_2(\gamma_2)\}
\end{aligned}\end{equation*}
be the fibre products over $\bar Q$ resp.\ $Q$, and, for $i\in\{1,2\}$, let
\begin{equation}\llabel{eq:defbarNi}
  \bar N_i=\pi_i^{-1}(Z)\leq\G_i.
\end{equation}
Suppose that $N_1$ is of type $\FP_k$ and that $N_2$ is of type
$\FP_l$.  
Then 
\begin{enumerate}[label=(\roman*)]
  \item \llabel{it:kernel-finiteness} the groups $\bar N_1$ and $\bar N_2$
  are of types $\FP_k$ and $\FP_l$, respectively,
  \item \llabel{it:novikov}
  $S(\bar P,P)\subseteq\Nov{k+l+1}{\bar P}$, and
  \item \llabel{it:finiteness} for every $d\in\N_0$ such that
  $\bar P$ is of type $\FP_d$ it holds that
  $P$ is of type $\FP_{\min\{d,k+l+1\}}$.
\end{enumerate}
\end{athm}
\begin{aproof}
First, observe that
\lref{eq:centr} and \lref{eq:defbarNi}
give short exact sequences
\[
\ses{
  1\to\bar N_1\to\G_1\xrightarrow{p\circ\pi_1}\bar Q\to1
} {
  1\to\bar N_2\to\G_2\xrightarrow{p\circ\pi_2}\bar Q\to1
}
\]
  Next, note that for all $(\gamma_1,\gamma_2)\in\bar P$ it holds that
  \[
    p(\pi_1(\gamma_1)\pi_2(\gamma_2)^{-1})=p\pi_1(\gamma_1)p\pi_2(\gamma_2)^{-1}=1,
  \]
  so $\pi_1(\gamma_1)\pi_2(\gamma_2)^{-1}\in\ker p=Z$.
  Let $\delta\colon\bar P\to Z$ satisfy for all $(\gamma_1,\gamma_2)\in\bar P$ that
\begin{equation}
  \llabel{eq:def_delta}
  \delta(\gamma_1,\gamma_2)=\pi_1(\gamma_1)\pi_2(\gamma_2)^{-1}.
\end{equation}
The assumption that $Z$ is central in $Q$ ensures
  for all $(\gamma_1,\gamma_2),(\gamma_1',\gamma_2')\in\bar P$ that
  \begin{equation*}
    \begin{aligned}
    \delta((\gamma_1,\gamma_2)(\gamma_1',\gamma_2'))
    &=
    \pi_1(\gamma_1\gamma_1')\pi_2(\gamma_2\gamma_2')^{-1}
    =
    \pi_1(\gamma_1)\pi_1(\gamma_1')\pi_2(\gamma_2')^{-1}\pi_2(\gamma_2)^{-1}
    \\&=\pi_1(\gamma_1)\delta(\gamma_1',\gamma_2')\pi_2(\gamma_2)^{-1}
    =\pi_1(\gamma_1)\pi_2(\gamma_2)^{-1}\delta(\gamma_1',\gamma_2')
    =\delta(\gamma_1,\gamma_2)\delta(\gamma_1',\gamma_2').
  \end{aligned}
  \end{equation*}
Hence, $\delta$ is a homomorphism and \lref{eq:def_delta} shows that
  $\ker\delta=P$.
For all $\gamma_1\in\bar N_1$, $\gamma_2\in\bar N_2$ one has 
$(\gamma_1,1),(1,\gamma_2)\in\bar P$ and
\[
  \delta(\gamma_1,1)=\pi_1(\gamma_1),
  \qquad
  \delta(1,\gamma_2)=\pi_2(\gamma_2)^{-1}.
\]
Identifying $\bar N_1$ with $\bar N_1\times 1$ and
$\bar N_2$ with $1\times\bar N_2$,
surjectivity of $\pi_1$ and $\pi_2$ shows that
\begin{equation}
  \delta(\bar N_1)=Z=\delta(\bar N_2),
  \qquad
  \ker(\delta|_{\bar N_1})=N_1,
  \qquad\text{and}\qquad
  \ker(\delta|_{\bar N_2})=N_2.
\end{equation}
\Cref{lem:abelian-quotient-fibre-product}, applied with $A=Z$, now establishes
\lref{it:kernel-finiteness,it:novikov,it:finiteness}.
\end{aproof}

\section{Iterated central extensions and virtually nilpotent quotients}
\label{sec:nilpotent}

We first iterate the central-extension step relative to an arbitrary base
quotient.  The nilpotent case is obtained
by taking the base quotient to be trivial; a finite-index argument then gives the
virtually nilpotent case.

\begin{athm}{theorem}{thm:nilpotent-fp}
Let $d,k,l\in\N_0$ satisfy $d\leq k+l+1$. Let
\[
\ses{
  	  1\to N_1\to \G_1\xrightarrow{\pi_1}Q\to1
} {
    1\to N_2\to \G_2\xrightarrow{\pi_2}Q\to1
}\]
be short exact sequences of groups, with $N_1$ of type $\FP_k$ and $N_2$ of
type $\FP_l$, let $c\in\N$ and for every
$j\in\{1,2,\ldots,c\}$ let
\begin{equation}\llabel{eq:central_ext}
  1\to Z_j\to Q_j\xrightarrow{p_j}Q_{j-1}\to1
\end{equation}
be a central extension with $Z_j$ finitely generated, assume that
$Q_c=Q$, let
\[
\begin{aligned}
	  \bar P&=
  \set{(\gamma_1,\gamma_2)\in\G_1\times\G_2\given
  (p_1\circ\dots\circ p_c\circ\pi_1)(\gamma_1)=(p_1\circ\dots\circ p_c\circ\pi_2)(\gamma_2)},\\
  P&=
  \set{(\gamma_1,\gamma_2)\in\G_1\times\G_2\given
  \pi_1(\gamma_1)=\pi_2(\gamma_2)}
\end{aligned}
\]
be the fibre products over $Q_0$ and $Q$, respectively, and assume that
$\bar P$ is of type $\FP_d$.
Then $P$ is of type $\FP_d$.
\end{athm}
\begin{aproof}
  For every $j\in\{0,1,\dots,c-1\}$ let
  $q_j=p_{j+1}\circ\dots\circ p_c$, let $q_c=\id_Q$,
  for every $j\in\{0,1,\dots,c\}$, $i\in\{1,2\}$ let
  $N_i^j = \ker(q_j\circ\pi_i)$ and
\[
\begin{aligned}
  P_j&=
  \set{(\gamma_1,\gamma_2)\in\G_1\times\G_2\given
  (q_j\circ\pi_1)(\gamma_1)=(q_j\circ\pi_2)(\gamma_2)}.
\end{aligned}
\]

For every $j\in\{0,1,\dots,c\}$, there are then short exact sequences
\begin{equation}\llabel{eq:seqs2}
\ses{
  	1\to N_1^j\to \G_1\xrightarrow{q_j\circ\pi_1}Q_j\to1
} {
    1\to N_2^j\to \G_2\xrightarrow{q_j\circ\pi_2}Q_j\to1
}
\end{equation}
Moreover, it holds for every $j\in\{1,2,\dots,c\}$ and $i\in\{1,2\}$ that
\begin{equation}\llabel{eq:kernels}
  N_i^{j-1}
  = \ker(q_{j-1}\circ\pi_i) 
  = \ker(p_{j}\circ q_{j}\circ\pi_i)
  = (q_{j}\circ\pi_i)^{-1}(\ker p_{j})
  = (q_{j}\circ\pi_i)^{-1}(Z_{j}).
\end{equation}
\Itref{prop:central-step}{it:kernel-finiteness}, applied to the sequences in
\lref{eq:seqs2} and the central extension in
\lref{eq:central_ext}, hence shows for every $j\in\{1,2,\dots,c\}$
such that
$N_1^{j}$ is of type $\FP_k$ and $N_2^{j}$ is of type $\FP_l$ that
$N_1^{j-1}$ is of type $\FP_k$ and $N_2^{j-1}$ is of type $\FP_l$.  Thus the
assumptions that $N_1^c=N_1$ is of type $\FP_k$ and $N_2^c=N_2$ is of type
$\FP_l$ and downward induction show that for every
$j\in\{0,1,\dots,c\}$ it holds that
$N_1^j$ is of type $\FP_k$ and $N_2^j$ is of type $\FP_l$.

\Itref{prop:central-step}{it:finiteness}
now establishes 
for every $j\in\{1,2,\dots,c\}$ such that  $P_{j-1}$ is of type $\FP_d$
that
$P_{j}$ is of type $\FP_d$.  Thus the assumption that $P_0=\bar P$ is of
type $\FP_d$ and induction show that for
every $j\in\{0,1,\dots,c\}$ it holds that $P_j$ is of type $\FP_d$.
In particular, $P_c=P$ is of type $\FP_d$.
\end{aproof}

A particularly important special case of the preceding theorem is
the case where $Q_0=1$, whence $Q$ is a nilpotent group.

\begin{athm}{corollary}{cor:nilpotent-fp}[Symmetric $n$-$(n+1)$-$(n+2)$ theorem
  for nilpotent quotients]
Let $d,k,l\in\N_0$ satisfy $d\leq k+l+1$, let
\[
\ses{
  	  1\to N_1\to \G_1\xrightarrow{\pi_1}Q\to1
} {
    1\to N_2\to \G_2\xrightarrow{\pi_2}Q\to1
}\]
be short exact sequences of groups.  Suppose that $\G_1$ and $\G_2$ are of type
$\FP_d$, that $N_1$ is of type $\FP_k$, that $N_2$ is of type $\FP_l$,
and that $Q$ is nilpotent.  Then the associated fibre product $P$ is of type
$\FP_d$.
\end{athm}
\begin{aproof}
  Assume without loss of generality that $d\geq 1$.
  Nilpotency of $Q$ means that there exists a $c\in\N$ and
  for every $j\in\{1,2,\dots,c\}$
  central extensions
  \begin{equation}
  1\to Z_j\to Q_j\xrightarrow{p_j}Q_{j-1}\to1
  \end{equation}
  with $Q_c=Q$ and $Q_0=1$.
  Since $\Gamma_1$ was assumed to be of type $\FP_d$ and 
  $d\geq 1$, it holds that $\Gamma_1$ is finitely generated and hence so
  are its successive nilpotent quotients $Q=Q_c,Q_{c-1},\dots,Q_1$.
  Thus, for all $j\in\{1,2,\dots,c\}$ the subgroups
  $Z_j$ are also finitely generated (see, e.g., \cref{it:nilpotent-fpinfty} in \cref{lem:finiteness-facts}).
  Since $Q_0=1$, the fibre product over $Q_0$ is equal to
  $\Gamma_1\times\Gamma_2$, which is of type $\FP_d$ by
  the assumption that $\Gamma_1$ and $\Gamma_2$ are of type $\FP_d$ and
\cref{it:ext-preserve-fp} in \cref{lem:finiteness-facts}.  The result now follows from
\cref{thm:nilpotent-fp}.
\end{aproof}

Passing to a finite-index nilpotent subgroup of the quotient gives the
corresponding result for virtually nilpotent quotients.

\begin{athm}{corollary}{cor:virtually-nilpotent-fp}[Symmetric $n$-$(n+1)$-$(n+2)$ theorem
  for virtually nilpotent quotients]
Let $d,k,l\in\N_0$ satisfy $d\leq k+l+1$, let
\[
\ses{
    1\to N_1\to \G_1\xrightarrow{\pi_1}Q\to1
} {
    1\to N_2\to \G_2\xrightarrow{\pi_2}Q\to1
}\]
be short exact sequences of groups.  Suppose that $\G_1$ and $\G_2$ are of type
$\FP_d$, that $N_1$ is of type $\FP_k$, that $N_2$ is of type $\FP_l$,
and that $Q$ is virtually nilpotent.  Then the associated fibre product $P$ is
of type $\FP_d$.
\end{athm}

\begin{aproof}
  Assume without loss of generality that $d\geq 1$.
Then $Q$ is finitely generated as a quotient of the
finitely generated group $\G_1$.  Choose a finite-index nilpotent subgroup
$Q^\circ\leq Q$ and for $i\in\{1,2\}$ put
\[
  \G_i^\circ=\pi_i^{-1}(Q^\circ).
\]
Then $\G_i^\circ$ has finite index in $\G_i$ and hence is of type $\FP_d$ by
\cref{it:finite-index} in \cref{lem:finiteness-facts}.  The restricted maps give short exact sequences
\[
\ses{
    1\to N_1\to \G_1^\circ\xrightarrow{\pi_1|}Q^\circ\to1
} {
    1\to N_2\to \G_2^\circ\xrightarrow{\pi_2|}Q^\circ\to1.
}
\]
Since $Q^\circ$ is nilpotent, \cref{cor:nilpotent-fp} shows that their
associated fibre product
\[
  P^\circ=
  \set{(\gamma_1,\gamma_2)\in\G_1^\circ\times\G_2^\circ\given
  \pi_1(\gamma_1)=\pi_2(\gamma_2)}
\]
is of type $\FP_d$.
The homomorphism
\[
  P\to Q,\   (\gamma_1,\gamma_2)\mapsto\pi_1(\gamma_1)=\pi_2(\gamma_2),
\]
is surjective, and $P^\circ$ is the preimage of $Q^\circ$.  Thus $P^\circ$
has finite index in $P$, so
applying \cref{it:finite-index} in  \cref{lem:finiteness-facts} again establishes that $P$ is
of type $\FP_d$.
\end{aproof}

\begin{athm}{corollary}{cor:virtually-nilpotent-fp-2}[Symmetric $n$-$(n+1)$-$(n+2)$ theorem
  for virtually nilpotent quotients]
Let $k,l,n,m\in\N_0$, let
\[
\ses{
    1\to N_1\to \G_1\xrightarrow{\pi_1}Q\to1
} {
    1\to N_2\to \G_2\xrightarrow{\pi_2}Q\to1
}\]
be short exact sequences of groups.  Suppose that $\G_1$ is of type $\FP_m$,
$\G_2$ is of type $\FP_n$, $N_1$ is of type $\FP_k$, and $N_2$ is of type $\FP_l$,
and assume that $Q$ is virtually nilpotent. Then the associated fibre product $P$ is
of type $\FP_{\min\{k+l+1,n,m\}}$.
\end{athm}
\begin{aproof}
  Let $d=\min\{k+l+1,n,m\}$. Then $d\leq k+l+1$ and $\G_1$ and $\G_2$ are of type
  $\FP_d$. \cref{cor:virtually-nilpotent-fp} now shows that $P$ is of type $\FP_d$.
\end{aproof}

Combining the homological result with the $1$-$2$-$3$ theorem and
Wall's criterion gives the homotopical version in \cref{cor:homotopical-virtually-nilpotent-fp} below.
\begin{athm}{theorem}{thm:123}[$1$-$2$-$3$ theorem; Bridson, Howie, Miller, \& Short~{\cite[Theorem~B]{BHMS2013}}]
Let
\[
\ses{
  1\to N_1\to \G_1\to Q\to1
} {
  1\to N_2\to \G_2\to Q\to1
}
\]
be short exact sequences of groups. Suppose that $N_1$ is finitely generated,
that $\G_1$ and $\G_2$ are finitely
presented, and that $Q$ is of type $\F_3$.
Then the associated fibre product is finitely presented.
\end{athm}

\begin{athm}{corollary}{cor:homotopical-virtually-nilpotent-fp}[Symmetric $n$-$(n+1)$-$(n+2)$ theorem for virtually nilpotent quotients, homotopical form]
Let $k,l,n,m\in\N_0$, let
\[
\ses{
    1\to N_1\to \G_1\xrightarrow{\pi_1}Q\to1
} {
    1\to N_2\to \G_2\xrightarrow{\pi_2}Q\to1
}\]
be short exact sequences of groups.  Suppose that $\G_1$
is of type $\F_m$, that $\G_2$ is of type $\F_n$,
that $N_1$ is of type $\FP_k$, that $N_2$ is of type $\FP_l$, and
that $Q$ is virtually nilpotent.  Then the associated fibre product $P$ is of
type $\F_{\min\{m,n,k+l+1\}}$.
\end{athm}
\begin{aproof}
Let $d=\min\{m,n,k+l+1\}$.
Since $\G_1$ and $\G_2$ are of type $\F_d$, they are of type
$\FP_d$.
\Cref{cor:virtually-nilpotent-fp-2}
therefore gives that $P$ is of type $\FP_d$.

If $d=0$, the conclusion is vacuous.
If $d=1$, then $P$ is finitely generated, hence of type
$\F_1$ (cf.\ \cref{lem:finiteness-facts}).
Suppose that $d\geq 2$. Then $k+l+1\geq 2$, so 
at least one of $k$ and $l$ is positive. Thus at least one
of $N_1$ and $N_2$ is finitely generated. 
Since $m\geq d$ and $n\geq d$, the groups $\G_1$ and $\G_2$ are
finitely presented. Their quotient $Q$ is then a finitely
generated virtually nilpotent group, hence of type $\F_\infty$
(see \cref{it:finite-index,it:nilpotent-fpinfty} in \cref{lem:finiteness-facts}).  
\Cref{thm:123} (possibly after interchanging the two short
exact sequences) thus establishes that $P$ is finitely presented.
The fact that $P$ is of type $\FP_d$ and Wall's criterion (see \cref{it:wall} in \cref{lem:finiteness-facts})
now ensure that $P$ is of type $\F_d$.
\end{aproof}

\section{Consequences and sharpness}

\label{sec:consequences}

Cohen and Shusterman have recently proved the Virtual Surjections Theorem below
in full generality, in both its homotopical and homological forms.
We record here that this theorem can be recovered from
the results proved above.

\begin{athm}{corollary}{cor:vsc}[Virtual Surjections Theorem, Cohen \& Shusterman {\cite[Theorem~1.3]{CohenShusterman2026}}]
  Let $m,n\in\N$ satisfy $2\leq m\leq n$, let
  $\Gamma_1,\Gamma_2,\dots,\Gamma_n$ be groups
  of type $\F_m$ (resp.\ $\FP_m$),
  let $\Gamma=\Gamma_1\times\Gamma_2\times\dots\times\Gamma_n$ be their direct product,
  let $G\leq \Gamma$
be a subgroup, and
assume that for all $i_1,i_2,\dots, i_m\in\{1,2,\dots,n\}$ 
with $i_1< i_2<\dots< i_m$ it holds that
the image of $G$ under the natural projection
$\Gamma\to\Gamma_{i_1}\times\dots\times\Gamma_{i_m}$ has finite index.
Then $G$ is of type $\F_m$ (resp.\ $\FP_m$).
\end{athm}
\begin{aproof}
  The fact that the homotopical statement involving the $\F_m$-properties follows from 
  \cref{cor:homotopical-virtually-nilpotent-fp} is
  \cite[Theorem~3.10]{Kuckuck2014}.
  In the homological setting the analogous implication was
  proved by Kochloukova and Lima~\cite[Theorem F]{KochloukovaLima2018}, so
  \cref{cor:virtually-nilpotent-fp-2} establishes the statement involving the $\FP_m$
  properties.
\end{aproof}

Next we show that the term $k+l+1$ in the conclusions of \cref{cor:virtually-nilpotent-fp-2,cor:homotopical-virtually-nilpotent-fp} cannot be improved, even when the quotient is infinite cyclic and the ambient groups are of type $\F_\infty$.

\begingroup
\begin{athm}{proposition}{prop:sharpness}
  Let $k,l\in\N_0$. Then there exist short exact sequences
  \begin{equation}\llabel{eq:ses}
    \ses{
    1\to N_1\to\Gamma_1\to \Z\to 1
    } {
    1\to N_2\to\Gamma_2\to \Z\to 1
    }
  \end{equation}
  such that $N_1$ is of type $\F_k$, $N_2$ is of type $\F_l$,
  $\Gamma_1$ and $\Gamma_2$ are of type $\F_\infty$,
  and the fibre product associated to \lref{eq:ses} is not of type $\FP_{k+l+2}$.
\end{athm}
\begin{aproof}
  Let $F$ be a free group on two generators
  and consider for every $s\in\N$
  the direct product $F^s=F\times F\times \dots\times F$
  with, for every $i\in\{1,2,\dots,s\}$, the $i$-th factor freely generated by
  $a_i,b_i\in F^s$ and the
  homomorphism $\phi_s\colon F^s\to\Z$ that satisfies for all
  $i\in\{1,2,\dots,s\}$ that $\phi_s(a_i)=1=\phi_s(b_i)$.
  The kernels of these homomorphisms are the Stallings--Bieri groups (see
  Bieri~\cite[Proposition~2.14]{Bieri1976}), which satisfy
  for all $s\in\N$ that
  \begin{equation*}
    \text{$\ker \phi_s$ is of type $\F_{s-1}$ but not of type $\FP_s$}.
  \end{equation*}
  Now let
  $P_{k,l}$ be the fibre product
  associated to the short exact sequences
  \begin{equation*}
    \ses{
      1\to \ker \phi_{k+1}\to F^{k+1}\xrightarrow{\phi_{k+1}} \Z\to1
    } {
      1\to\ker(-\phi_{l+1})\to F^{l+1}\xrightarrow{-\phi_{l+1}} \Z\to1
    }
  \end{equation*}
  Then
  \begin{equation*}
    P_{k,l} = \set{(\gamma_1,\gamma_2)\in F^{k+1}\times F^{l+1}\given \phi_{k+1}(\gamma_1)+ \phi_{l+1}(\gamma_2)=0}.
  \end{equation*}
  Therefore, $P_{k,l}$ is simply the kernel of $\phi_{k+l+2}$ (under
  the identification of $F^{k+1}\times F^{l+1}$ with $F^{k+l+2}$)
  and hence not of type $\FP_{k+l+2}$.
\end{aproof}
\endgroup

\section{Fibre products over quotients with coherent integral group ring}
\label{sec:coherent-quotients}

In this section, we will prove a weak version of 
\cref{cor:virtually-nilpotent-fp-2}
for fibre products over quotients whose integral group rings are coherent, which includes
virtually polycyclic groups as special cases. In this section, all modules will be left
modules.

\begin{definition}
  A ring $R$ is called (left) \emph{coherent} if every finitely generated (left) ideal of $R$ is finitely presented.
\end{definition}

We will use the following elementary 
reformulations.
\begin{athm}{lemma}{lem:coherent-reformulations}
  Let $R$ be a ring. Then the following are equivalent:
  \begin{enumerate}[(i)]
    \item \llabel{it:coherent} $R$ is coherent.
    \item \llabel{it:fpinfty} Every finitely presented $R$-module is of type $\FP_\infty$.
    \item \llabel{it:fpkernel} Every homomorphism between finitely presented $R$-modules
    has finitely presented kernel.
    \item \llabel{it:fgkernel} Every homomorphism between finitely generated free $R$-modules
    has finitely generated kernel.
  \end{enumerate}
\end{athm}
\begin{aproof}
  The equivalence \lref*{it:coherent}$\iff$\lref*{it:fpinfty} and the
  implication \lref*{it:fpinfty}$\implies$\lref*{it:fpkernel}
  are proved, e.g., in Bravo, Gillespie, and Hovey~\cite[Proposition~2.1 and following discussion]{bravo2014stable}.
  The implication \lref*{it:fpkernel}$\implies$\lref*{it:fgkernel} is immediate.
  Finally, if \lref{it:fgkernel} holds, then for every finitely presented $M$
  a finite presentation $R^m\to R^n\to M\to 0$ with $m,n\in\N$ can be
  successively extended to the left to a resolution of $M$ by finitely generated free $R$-modules.
  This shows \lref*{it:fgkernel}$\implies$\lref*{it:fpinfty}.
\end{aproof}

The following gives a collection of some classes of groups 
whose integral group rings are known to be coherent.

\begin{athm}{proposition}{prop:coherent-group-rings} Let $G$, $G_1$, $G_2$, $H$ be groups. Then the following hold:
  \begin{enumerate}[(i)]
    \item \llabel{it:polycyclic} If $G$ is virtually polycyclic, then $\Z G$ is Noetherian and hence coherent.
    \item \llabel{it:amalgamated} If $\Z [G_1]$ and $\Z [G_2]$ are coherent
    and $\Z H$ is Noetherian, then the integral group ring 
    of every amalgamated free product $G_1*_H G_2$ and every HNN extension $G_1*_H$ is coherent.
    \item \llabel{it:free} If $G$ is a finitely generated free or surface group, then $\Z G$ is coherent.
  \end{enumerate}
\end{athm}
\begin{aproof}
  \Lref{it:polycyclic} is due to Hall~\cite{Hall1954}.
  \Lref{it:amalgamated} is essentially due to Waldhausen~\cite[Theorem~19.1]{Waldhausen1978},
  see also {\AA}berg~\cite{Aberg1982} and Lam~\cite{Lam1977}.
  \Lref{it:free} follows from \lref{it:polycyclic} and \lref{it:amalgamated}
  as all finitely generated free groups can, e.g., be obtained by starting with the trivial group 
  and forming iterated HNN extensions with trivial edge group
  and every infinite, non-free
  surface group can be obtained as an HNN extension of a free group
  with cyclic edge group.
\end{aproof}

\begin{definition}\label{def:lhs-controlled}
  Given a short exact sequence
\begin{equation}\label{eq:lhs-ses}
  1\to N\to \Gamma\to Q\to1
\end{equation}
and $r\in\N$, we say that $N$ is \emph{$r$-LHS-controlled over $Q$} if, with
the quotient (left) $Q$-action on $H_*(N;\Z)$,
\begin{enumerate}[(i)]
  \item\label{it:lhs-controlled-1} it holds for all $q\in\{0,1,\dots,r-1\}$ that $H_q(N;\Z)$ is of type $\FP_\infty$ over $\Z Q$ and
  \item\label{it:lhs-controlled-2} it holds that $H_r(N;\Z)_Q$ is finitely generated as an abelian group.
\end{enumerate}
\end{definition} 

As indicated by the name, this definition is
rather transparently chosen to ensure that the terms 
relevant to the computation of $H_*(\Gamma;\Z)$ up to a certain degree
on the second page of the 
Lyndon--Hochschild--Serre spectral sequence associated to
\cref{eq:lhs-ses} 
are finitely generated. So the following is almost immediate.

\begin{athm}{lemma}{lem:lhs-transfer}
Let $r\in\N$, let
\[
  1\to N\to\Gamma\to Q\to1
\]
be a short exact sequence and assume that $N$ is $r$-LHS-controlled over $Q$.
Then
\begin{equation}
  \llabel{eq:conc}
  \forall j\in\{0,1,\dots,r\}\colon \text{$H_j(\Gamma;\Z)$ is
  finitely generated}.
\end{equation}
\end{athm}
\begin{aproof}
Consider the Lyndon--Hochschild--Serre spectral sequence
\[
  E^2_{p,q}=H_p(Q;H_q(N;\Z))
  \Longrightarrow H_{p+q}(\Gamma;\Z).
\]
The assumption that $N$ is $r$-LHS-controlled over $Q$
implies that for every $q\in\{0,1,\dots,r-1\}$ it holds that $H_q(N;\Z)$ has a resolution
by finitely generated projective $\Z Q$-modules.
Tensoring such a resolution with the trivial right $\Z Q$-module
$\Z$ give a chain complex of finitely generated abelian groups.
Hence, it holds for every $p,q\in\N_0$ with $q<r$ that 
the homology groups $E^2_{p,q}=H_p(Q;H_q(N;\Z))$
are finitely generated.
Furthermore, by \cref{it:lhs-controlled-2}
in \cref{def:lhs-controlled}
the assumption that $N$ is $r$-LHS-controlled over $Q$
implies that $E^2_{0,r}=H_r(N;\Z)_Q$ is finitely generated.
Since for every $p,q,s\in\N_0$
with $s\geq 3$ the term $E^{s}_{p,q}$ is a subquotient of
$E^{s-1}_{p,q}$, it follows that likewise the terms of
total degree at most $r$ on the $E^\infty$ page are
finitely generated. This proves \lref{eq:conc}.
\end{aproof}
Less obviously, when $\Z Q$ is coherent,
there is a kind of converse to \cref{lem:lhs-transfer}.

\begin{athm}{lemma}{lem:coherent-lhs-transfer}
Let $r\in\N$, let
\[
  1\to N\to\Gamma\xrightarrow\pi Q\to1
\]
be a short exact sequence, assume that $\Gamma$ is of type $\FP_r$,
  and assume that $\Z Q$ is coherent. Then $N$ is $r$-LHS-controlled over $Q$.
\end{athm}
\begin{aproof}
Choose a free
resolution $F_*\to\Z$ of the trivial left $\Z\Gamma$-module $\Z$ such that $F_q$
is finitely generated for $q\leq r$. Regard $\Z Q$ as a
$(\Z Q,\Z\Gamma)$-bimodule via the quotient map $\pi$ and put
\[
  (C_*,d_*)=\Z Q\tensor_{\Z\Gamma}F_*.
\]
By Shapiro's Lemma, there is a chain isomorphism
$\Theta\colon C_*=\Z Q\tensor_{\Z\Gamma}F_*\to \Z\tensor_{\Z N}F_*$, 
such that for every $\gamma\in\Gamma$, $i\in\N_0$, $f\in F_i$
it holds that $\Theta(\pi(\gamma)\tensor f)=1\tensor \gamma f$
and $\Theta$ is $Q$-equivariant (with $Q$ acting on the left factor of $\Z Q\tensor_{\Z\Gamma} F_*$
and on the right factor of $\Z\tensor_{\Z N}F_*$ via lifting to $\Gamma$).
This map therefore induces for every $q\in\N_0$ a $Q$-equivariant isomorphism on homology groups
\begin{equation}\llabel{eq:shapiro}
  H_q(C_*) \cong H_q(\Z\tensor_{\Z N}F_*) \cong H_q(N;\Z).
\end{equation}
Moreover, it holds for all $q\leq r$ that
$C_q$ is a finitely generated free, and therefore finitely
presented, $\Z Q$-module. The assumption that
$\Z Q$ is coherent and \cref{lem:coherent-reformulations}
hence ensure for all $q<r$ that
$\ker d_q$ and $\im d_{q+1}$ are likewise finitely presented.
Combining this with \lref{eq:shapiro} establishes that
\begin{equation}\llabel{eq:finpres}
  \forall q\in\{0,1,\dots,r-1\}\colon \text{$H_q(N;\Z)$ is finitely presented as a $\Z Q$-module}.
\end{equation}
Coherence of $\Z Q$ and \cref{lem:coherent-reformulations} therefore imply that 
\begin{equation}\llabel{eq:FPinfty}
  \forall q\in\{0,1,\dots,r-1\}\colon \text{$H_q(N;\Z)$ is of type $\FP_\infty$ as a $\Z Q$-module}.
\end{equation}
Moreover, coherence of $\mathbb ZQ$ 
and \cref{lem:coherent-reformulations} ensure that $\ker(d_r\colon C_r\to C_{r-1})$ is finitely generated over $\mathbb ZQ$, since $C_r$ and $C_{r-1}$ are finitely generated free. Hence
\[
H_r(N;\mathbb Z)\cong\ker d_r/\operatorname{im}d_{r+1}
\]
is finitely generated over $\mathbb ZQ$ and so $H_r(N;\mathbb Z)_Q$ is finitely generated over $\mathbb Z$. 
This and \lref{eq:FPinfty} establish that $N$ is
$r$-LHS-controlled over $Q$.
\end{aproof}

We will use the two preceding lemmas to transfer finiteness assumptions on the
groups in the two short exact sequences
\begin{equation}\label{eq:ses-coherent}
  \ses{
    1\to N_1\to\Gamma_1\to Q\to1
  }{
    1 \to N_2\to\Gamma_2\to Q\to1
  }
\end{equation}
to the short exact sequence
\begin{equation*}
  1\to N_1\times N_2\to P\to Q\to 1
\end{equation*}
for the fibre product $P$ associated to \cref{eq:ses-coherent}.
This will require analyzing the terms in the Künneth 
sequence for computing the homology of $N_1\times N_2$.
The next two lemmas first supply the necessary finiteness properties of these terms.

\begin{athm}{lemma}{lem:tensor-finpres}
Let $Q$ be a group, let $A$ be a left $\Z Q$-module 
that is finitely
generated as an abelian group, and let $M$ be a finitely presented
left $\Z Q$-module. Then the modules $A\tensor_\Z M$ and $M\tensor_\Z A$ with the diagonal $Q$-action
are finitely presented over $\Z Q$.
\end{athm}

\begin{aproof}
The assumption that $M$ is finitely presented means that there exist
$u,v\in\N$ and an exact sequence of left $\Z Q$-modules
\begin{equation*}
  (\Z Q)^u\to (\Z Q)^v\to M\to0.
\end{equation*}
Applying $A\tensor_\Z -$ yields an exact sequence
\begin{equation}\llabel{eq:presA}
  (A\tensor_\Z \Z Q)^u\to (A\tensor_\Z \Z Q)^v\to A\tensor_\Z M\to0.
\end{equation}
Here, $Q$ acts diagonally on $A\tensor_\Z \Z Q$.
This module is isomorphic to the left $\Z Q$-module
$\Z Q\tensor_\Z A$ where $Q$ acts only on the left factor.
Indeed, the
map $A\times Q\to \Z Q\tensor_\Z A$, $(a, q)\mapsto q\tensor q^{-1}a$
induces a homomorphism $\Phi\colon A\tensor_\Z \Z Q\to \Z Q\tensor_\Z A$ of left $\Z Q$-modules
which has an inverse induced by the map $Q\times A\to A\tensor_\Z \Z Q$,
$(q,a)\mapsto qa\tensor q$.
Since $A$ is finitely generated as an abelian group,
$\Z Q\tensor_{\Z}A$ is finitely presented over $\Z Q$
and, via the isomorphism $\Phi$, so is
$A\tensor_\Z \Z Q$.
Therefore, in particular,
the first module in \lref{eq:presA} is finitely generated and
the second is finitely presented. Thus
$A\tensor_\Z M$ is finitely presented. Since
$M\tensor_\Z A\cong A\tensor_\Z M$ as left $\Z Q$-modules, the conclusion follows.
\end{aproof}

\begin{athm}{lemma}{lem:tor-finpres}
Let $Q$ be a group such that $\Z Q$ is coherent, let $A$ be a left $\Z Q$-module 
that is finitely
generated as an abelian group, and let $M$ be a finitely presented
left $\Z Q$-module. Then the modules $\Tor^{\Z}_1(A,M)$ and $\Tor^{\Z}_1(M,A)$ with the diagonal $Q$-action
are finitely presented over $\Z Q$.
\end{athm}

\begin{aproof}
Let $T$ be the torsion subgroup of $A$, which is characteristic in $A$
and hence $Q$-invariant and furthermore finite by the 
assumption that $A$ is finitely generated as an abelian group.
The long exact $\Tor^{\mathbb Z}$-sequence associated to
$ 0\to T\to A\to A/T\to0 $
gives a natural $Q$-equivariant isomorphism
\begin{equation}\llabel{eq:tor-torsion}
  \Tor^{\mathbb Z}_1(T,M) \cong \Tor^{\mathbb Z}_1(A,M).
\end{equation}
Let $\Z^{(T)}$ be the free abelian group on the basis
$(e_t)_{t\in T}$, with
the $Q$-action induced by permutation of these basis elements
and let $\pi\colon \Z^{(T)}\to T$ be the homomorphism that satisfies
for every $t\in T$ that $\pi(e_t)=t$.
Then
\[
  0\to \ker\pi\xrightarrow\iota \Z^{(T)}\xrightarrow\pi T\to0
\]
is a $Q$-equivariant short exact sequence with $\Z^{(T)}$ and
$\ker\pi$ free abelian groups.
It follows that
\begin{equation}\llabel{eq:ker}
  \Tor^{\Z}_1(T,M)
  =\ker\bigl(
    (\ker\pi)\tensor_{\Z}M\xrightarrow{\iota\tensor_\Z M} \Z^{(T)}\tensor_{\Z}M
  \bigr).
\end{equation}
Both $\Z^{(T)}$ and $\ker\pi$ are finitely generated as abelian
groups, so \cref{lem:tensor-finpres}
shows that $(\ker\pi)\tensor_{\Z}M$ and $\Z^{(T)}\tensor_{\Z}M$
are finitely presented as $\Z Q$-modules.
Therefore, \lref{eq:ker}
and the assumption that $\Z Q$ is coherent
prove that $\Tor^{\Z}_1(T,M)$ is finitely presented over $\Z Q$.
Using \lref{eq:tor-torsion}
and the fact that
$\Tor^{\Z}_1(M,A)\cong \Tor^{\Z}_1(A,M)$ as left $\Z Q$-modules,
the conclusion follows.
\end{aproof}

\begin{athm}{lemma}{lem:coherent-lhs-product}
Let $r\in\N$, $k_1,k_2\in\N_0$ satisfy
$
  r\leq k_1+k_2+1
$.
Let
\[
\ses{
    1\to N_1\to \G_1\xrightarrow{\pi_1}Q\to1
} {
    1\to N_2\to \G_2\xrightarrow{\pi_2}Q\to1
}
\]
be short exact sequences, let $P$ be the associated fibre product,
assume that $\Z Q$ is coherent, assume for 
$i\in\{1,2\}$ that 
\begin{equation}\llabel{eq:assNi}
  \forall j\in\{0,1,\dots,k_i\}\colon \text{$H_j(N_i;\Z)$ is finitely generated over $\Z$},
\end{equation}
and assume that $N_1$ and $N_2$ are
$r$-LHS-controlled over $Q$.
Then there is a short exact sequence
\[
  1\to  N_1\times N_2\to  P\to  Q\to 1
\]
and  $N_1\times N_2$ is $r$-LHS-controlled over $Q$.
\end{athm}

\begin{aproof}
Throughout this proof, let $K=N_1\times N_2$.
There is a short exact Künneth sequence
(see, e.g., Brown~\cite[Proposition~I.0.8 and Corollary~V.5.8]{Brown1982})
\begin{equation}\llabel{eq:kuenneth}
\begin{aligned}
0\to 
\bigoplus_{a+b=j}
H_a(N_1;\Z)\tensor_{\Z}H_b(N_2;\Z)
\to  H_j(K;\Z)
\to  \!\!
\bigoplus_{a+b=j-1}\!\!
\Tor^{\Z}_1(H_a(N_1;\Z),H_b(N_2;\Z))
\to 0
\end{aligned}
\end{equation}
which is natural in both factors and therefore $Q$-equivariant with each
summand in the outer terms carrying the diagonal $Q$-action.
Note that the assumption that 
$N_1$ and $N_2$ are $r$-LHS-controlled over $Q$
implies for $i\in\{1,2\}$ that 
\begin{equation}\llabel{eq:fpinfty}
  \forall j\in\{0,1,\dots,r-1\}\colon \text{$H_j(N_i;\Z)$ is of type $\FP_\infty$ over $\Z Q$}.
\end{equation}
Moreover, \lref{eq:assNi} and the assumption that $r\leq k_1+k_2+1$
ensure that
\begin{equation}
  \begin{multlined}
  \forall a,b\in\N_0\colon a+b\leq r \implies (\text{$H_a(N_1;\Z)$ or $H_b(N_2;\Z)$ is finitely generated}\\
  \text{as an abelian group}).
  \end{multlined}
\end{equation}
Combining this, \lref{eq:fpinfty}, and the assumption that $\Z Q$ is coherent
with \cref{lem:tensor-finpres,lem:tor-finpres} ensures for all 
$a,b\in\N_0$ with $a+b\leq r$ and $\max\{a,b\}<r$ that
\begin{multline}\llabel{eq:finite}
  \text{$H_a(N_1;\Z)\tensor_{\Z}H_b(N_2;\Z)$ and $\Tor^{\Z}_1(H_a(N_1;\Z),H_b(N_2;\Z))$
  are finitely presented}\\\text{ as $\Z Q$-modules}.
\end{multline}
Coherence of $\Z Q$ and \lref{eq:kuenneth}
then show that
\begin{equation}\llabel{eq:lhs1}
  \forall j\in\{0,1,\dots,r-1\}\colon \text{$H_j(K;\Z)$ is of type $\FP_\infty$ over $\Z Q$}.
\end{equation}
Since the coinvariants functor $(-)_Q$ is right exact and commutes with
direct sums, \lref{eq:kuenneth} shows that there is an exact sequence
 \begin{equation}\llabel{eq:kuenneth2}
\begin{aligned}
\bigoplus_{a+b=r}
\br[\big]{H_a(N_1;\Z)\tensor_{\Z}H_b(N_2;\Z)}_Q
\to  H_r(K;\Z)_Q
\to  \!\!
\bigoplus_{a+b=r-1}\!\!
\br[\big]{\Tor^{\Z}_1(H_a(N_1;\Z),H_b(N_2;\Z))}_Q
\to 0.
\end{aligned}
\end{equation}
Furthermore, there are $Q$-equivariant isomorphisms
\begin{equation}
  \begin{aligned}
    &H_r(N_1;\Z)\tensor_{\Z}H_0(N_2;\Z)\cong H_r(N_1;\Z)\\
    \text{and}\qquad &H_0(N_1;\Z)\tensor_{\Z}H_r(N_2;\Z)  \cong H_r(N_2;\Z).\qquad\qquad
  \end{aligned}
\end{equation}
It follows that
\begin{multline*}
\bigoplus_{a+b=r} \br[\big]{H_a(N_1;\Z)\tensor_{\Z}H_b(N_2;\Z)}_Q\\[-0.5em]
\cong
H_r(N_2;\Z)_Q\oplus\pr[\Bigg]{\bigoplus_{\substack{a+b=r\\\max\{a,b\}<r}} \br[\big]{H_a(N_1;\Z)\tensor_{\Z}H_b(N_2;\Z)}_Q}\oplus H_r(N_1;\Z)_Q.
\end{multline*}
The assumption that $N_1$ and $N_2$ are $r$-LHS-controlled over $Q$ and
\lref{eq:finite} therefore show that 
\begin{equation}\llabel{eq:tensors}
\text{$\bigoplus_{a+b=r}
\br[\big]{H_a(N_1;\Z)\tensor_{\Z}H_b(N_2;\Z)}_Q$ is finitely generated over $\Z$}.
\end{equation}
Furthermore, \lref{eq:finite} ensures that
\begin{equation*}
  \text{$
  \bigoplus_{a+b=r-1}\!\!
\br[\big]{\Tor^{\Z}_1(H_a(N_1;\Z),H_b(N_2;\Z))}_Q$
is finitely generated over $\Z$}.
\end{equation*}
This, \lref{eq:tensors}, and \lref{eq:kuenneth2}
establish that
\[
  \text{$H_r(K;\Z)_Q$ is finitely generated over $\Z$.}
\]
This and \lref{eq:lhs1}
prove that $K$ is $r$-LHS-controlled over $Q$.
\end{aproof}

We can now state the coherent-quotient theorem in a form that is stronger than
the one obtained by assuming finiteness properties of the kernels.

\begin{athm}{theorem}{thm:coherent-quotient-homology}
Let $k,l,m,n\in\N_0$, let
\begin{equation}\llabel{eq:ses}
\ses{
    1\to N_1\to \G_1\xrightarrow{\pi_1}Q\to1
} {
    1\to N_2\to \G_2\xrightarrow{\pi_2}Q\to1
}
\end{equation}
be short exact sequences of groups, let $P$ be the associated fibre
product, assume that $\G_1$ is of type $\FP_m$ and that $\G_2$ is of type $\FP_n$,
assume that
\[
\begin{aligned}
  &\forall j\in\{0,1,\dots,k\}\colon \text{$H_j(N_1;\Z)$ is finitely generated}\\
  \text{and}\qquad&\forall j\in\{0,1,\dots,l\}\colon \text{$H_j(N_2;\Z)$ is finitely generated},\qquad
\end{aligned}
\]
and assume that $\Z Q$ is coherent.
Then 
\begin{equation}\llabel{eq:conc}
 \forall j\in\{0,1,\dots,\min\{m,n,k+l+1\}\}\colon
 \text{$H_j(P;\Z)$ is finitely
generated}.
\end{equation}
\end{athm}

\begin{aproof}
  Throughout this proof let $r=\min\{m,n,k+l+1\}$
  and assume without loss of generality that $r\geq1$.  Since
$r\leq \min\{m,n\}$, \cref{lem:coherent-lhs-transfer} applied to each of the two
short exact sequences in \lref{eq:ses} shows that $N_1$ and $N_2$ are $r$-LHS-controlled over
$Q$.  Since $r\leq k+l+1$, \cref{lem:coherent-lhs-product} now shows that
$K=N_1\times N_2$ is $r$-LHS-controlled over $Q$.  Applying
\cref{lem:lhs-transfer} to the short exact sequence
\[
  1\to  K\to  P\to  Q\to 1,
\]
yields the conclusion \lref{eq:conc}.
\end{aproof}

\subsubsection*{Acknowledgements}
This work has been supported by the Ministry of Culture and Science NRW as part of the
Lamarr Fellow Network. 

The author thanks ChatGPT (GPT-5.5 Pro \& GPT-5.6 Pro) for extensive mathematical discussions during the development of this work. These discussions contributed substantially to the discovery and refinement of several arguments, the exploration of proof strategies, and the identification of gaps. Final formulations,  exposition, and wording are entirely the author's and all mathematical statements and their proofs are the sole responsibility of the author.

\bibliographystyle{acm}
\bibliography{fibre_products_nilpotent}

\end{document}